\documentclass{article}

\usepackage[english]{babel}
\usepackage{graphicx}
\usepackage{subcaption}
\usepackage{amsmath}
\DeclareMathOperator*{\esssup}{ess\,sup}
\usepackage{amsthm}
\usepackage{amssymb}
\newtheorem{theorem}{Theorem}[section]

\newtheorem{proposition}{Proposition}[section]
\newtheorem{corollary}[theorem]{Corollary}
\theoremstyle{definition}
\newtheorem{definition}[theorem]{Definition}

\theoremstyle{definition}

\theoremstyle{remark}

\usepackage{float}
\usepackage{enumitem}
\usepackage{accents}

\usepackage[bbsymbols]{dsserif}
\newcommand{\lPi}{\rotatebox[origin=c]{180}{$\bbPi$}}
\newcommand{\uPi}{\bbPi}
\usepackage{natbib}
\usepackage[letterpaper,top=2cm,bottom=2cm,left=3cm,right=3cm,marginparwidth=1.75cm]{geometry}

\usepackage{amsmath}
\usepackage{graphicx}
\usepackage[colorlinks=true, allcolors=blue]{hyperref}

\title{Efficient inferential models and their fiducial connections}
\author{
    Max Raner\\[4pt]
    \small Department of Mathematics, Uppsala University, Sweden\\
    \small \texttt{max.raner@math.uu.se}
}

\begin{document}
\maketitle

\begin{abstract}
Valid inferential models (IMs) provide frequentist calibration guarantees, but validity alone does not determine which construction is most informative. This paper studies finite-sample efficiency through concentration risk, the expected volume under a possibility contour with respect to a specified parameter-volume measure. For a smooth, fully identified association, a weighted Neyman--Pearson argument identifies the generalized fiducial density as the optimal score for minimizing parameter-space cut volume subject to fixed-data auxiliary compatibility. The resulting contour is the density-based maximum-specificity transform of the fiducial distribution, giving an explicit volume interpretation of the fiducial efficiency ceiling. Its auxiliary ordering can depend on the observation, however, so this fixed-data optimum need not be sampling valid.

We develop two approaches to constructing valid IMs guided by this efficiency criterion. The first minimizes design-averaged concentration risk among single, data-independent auxiliary-compatible contours. The design measure specifies where efficiency is prioritized, without requiring a prior interpretation or restricting the validity guarantee. The second localizes the parameter-volume score to candidate-specific auxiliary contours, yielding a normalized, valid IM equivalent to sampling validification of the relative fiducial likelihood. A parallel sample-space construction recovers the ordinary likelihood-ratio IM, clarifying the geometric distinction between the two. Direct validification of the ceiling contour provides a complementary repair.

The fiducial connection extends beyond the optimizing construction: exact, data-independent predictive contours share the same fiducial distribution as a saturating inner probability, while producing different geometries of calibrated confidence regions. Design averaging therefore selects an efficient member of this family without changing the underlying fiducial law. Examples in a two-component measurement model and an exponential scale model illustrate the efficiency consequences of contour choice. In particular, a sampling-valid construction can improve on the raw ceiling by leaving the fixed-data auxiliary-compatible class over which it is optimal.
\end{abstract}

\noindent\textbf{Keywords:} Inferential models, generalized fiducial inference,
possibility theory, finite-sample efficiency, confidence sets.

\section{Introduction}

Valid inferential models (IMs) provide prior-free, posterior-like uncertainty summaries about unknown parameters in the form of lower and upper probabilities. Validity gives these assignments an assertion-wise frequentist calibration guarantee \citep{Martin2013}, controlling how often a true assertion receives small upper probability or, equivalently, a false assertion receives large lower support. The same summary thus yields tests with controlled type~I error and confidence regions with at least nominal coverage. Validity alone, however, leaves the choice of IM open. Just as several confidence-set procedures may be available for a given problem, different IMs can satisfy the same calibration requirement while differing substantially in informativeness. Of central importance, then, is the construction of statistically efficient IMs for the problem at hand.

In the original IM formulation, a data-generating association connects the observations and the unknown parameter through an auxiliary variable with a known distribution. Uncertainty about the unobserved auxiliary value is represented by a compatible predictive random set in auxiliary space and propagated through the association to obtain uncertainty about the parameter \citep{martin2015inferential}. Admissibility arguments motivate restricting attention to nested predictive random sets, whose induced upper probabilities are possibility measures \citep{martinhandbook,martin2015inferential}. This permits an equivalent formulation in terms of possibility contours (akin to a probability density), but leaves open the choice of association and, for a given association, the choice of predictive contour. More recent IM constructions obtain calibrated possibility contours directly from the sampling model, notably through the relative likelihood, without explicitly selecting a predictive random set \citep{Martin2025}.

Efficiency can be investigated from several complementary perspectives. \citet{Martin2025BvM} establish the asymptotic efficiency of the relative-likelihood IM under the regularity conditions of their possibilistic Bernstein--von Mises theorem. The present paper concerns the choice of an IM under a specified measure of finite-sample efficiency. We consider concentration risk, a particular case of the Choquet risk proposed by \citet{raner2026efficiency}. For a possibility contour $\pi_x$ and a parameter-volume measure $\nu$, the corresponding concentration loss is the volume under the contour, $L(\pi_x)=\int_\Theta \pi_x(\vartheta)\,\nu(d\vartheta)$. Concentration risk is its repeated-sampling expectation. Since the cuts $C_\alpha(x)=\{\vartheta:\pi_x(\vartheta)\geq\alpha\}$ of a valid contour are confidence regions \citep{Martin2021}, the risk has a classical interpretation as expected confidence-region volume integrated over all levels.

In this paper we return to the original IM construction and use this concentration criterion to guide the choice of predictive contour for a given association. When an auxiliary density is available, \citet{martinhandbook} recommend the contour supplied by the density-based \emph{maximum-specificity principle} \citep{DuboisEtAl2004}. This construction has a Neyman--Pearson interpretation as minimizing auxiliary volume subject to retaining a prescribed auxiliary probability. A small auxiliary set, however, can have a large image in parameter space. The concentration criterion suggests minimizing the volume of this image directly.

Recent IM developments have renewed interest in connections among Bayesian, fiducial, and frequentist inference \citep{reIM,hannig_demyst2025}. Here the IM--fiducial connection arises from the efficiency problem itself. For a smooth, fully identified association after any relevant dimension reduction, we first minimize parameter-space contour-cut volume subject to fixed-data auxiliary compatibility. This gives a weighted Neyman--Pearson set-selection problem \citep{lehmann2005testing}. We show that its optimal score is precisely the density, relative to $\nu$, of the generalized fiducial distribution induced by the inverse association \citep{hannig2009generalized}. The resulting contour is the density-based maximum-specificity transform of this distribution. Its cuts are fiducial highest-density regions of minimum $\nu$-volume at their prescribed content. \citet{martinhandbook} and \citet{hannig_demyst2025} relate fiducial probability to bounds on IM belief. Here the ceiling has an explicit volume interpretation within the class of fixed-data auxiliary-compatible contours.

Fixed-data auxiliary compatibility, however, need not imply sampling validity. When the optimal ordering depends on the observation, the auxiliary draw that generates the observation also determines the ordering by which that draw is ranked. The fiducial ceiling is consequently an efficiency benchmark within the fixed-data auxiliary compatible class of contours, but it need not itself define a valid IM. This distinction motivates the central constructive question of the paper: how can validity be restored to the ceiling contour and fiducial ranking, to obtain valid and efficient inferential models? 

We pursue two approaches. The first uses a single auxiliary-compatible contour chosen independently of the realized data. Averaging concentration risk over a design measure $w$ on the parameter space gives the optimal contour within this class. Sampling validity then follows from auxiliary compatibility. The design measure specifies where efficiency is to be prioritized and need not represent prior beliefs about the true parameter. In particular, efficiency can be targeted towards a scientifically relevant region while validity holds throughout the parameter space.

The second approach replaces a single common predictive contour by a candidate-indexed family of contours obtained through what we call \emph{localization}. For each candidate parameter value, we follow the observations it could generate and evaluate the ceiling's parameter-volume score along these hypothetical datasets. Normalizing and calibrating the local scores yields candidate-specific predictive contours, which are combined through the observed association into a normalized, valid IM. Equivalently, this random-set construction validifies the \emph{relative fiducial likelihood}, the fiducial density divided by its maximum at the given observation. A parallel localization based on sample-space volume recovers the ordinary validified likelihood-ratio IM. We also consider direct validification of the ceiling contour as a complementary repair.

These constructions provide an efficiency-oriented perspective on the saturating inner-probability interpretation of \citet{martin2023fiducial,reIM} and the fiducial confidence-curve interpretation of \citet{hannig_demyst2025}. For a fixed, fully identified association, exact, data-independent auxiliary contours generate sampling-calibrated IMs that share the same fiducial distribution as a saturating inner probability. Their nested confidence regions can differ in shape, yet each region has fiducial content equal to its nominal coverage. The design-averaged optimization selects an efficient member of this class without changing the underlying fiducial law.

We illustrate the constructions in a two-component measurement model with additive and multiplicative error \citep{Rocke} and in an exponential scale model. The measurement example shows that a sampling-valid repair can improve on the concentration of the raw ceiling, since sampling validity permits contours outside the fixed-data auxiliary-compatible class over which the ceiling is optimal. In the scale example, the optimal auxiliary ordering is data-independent, so the ceiling is already exactly calibrated and coincides with the proposed repairs. Auxiliary maximum specificity and validified relative likelihood also give exact coverage at every level. For $n=5$ and ordinary parameter length, their concentration risks exceed that of the ceiling by about $60\%$ and $8.7\%$, respectively, demonstrating that the choice of contour—even among reasonable alternatives—can have large efficiency consequences.

The remainder of the paper is organized as follows. Section~\ref{sec:background} reviews the theory on predictive contour constructions, auxiliary compatibility, and validity. Section~\ref{sec:auxcontours} derives the fixed-data efficiency ceiling and the optimal design-averaged construction. Section~\ref{sec:samplingrepair} develops the validity-repairs and their relationship with the validified likelihood-ratio IM. Section~\ref{sec:connections} establishes the fiducial representations and further literature connections. Section~\ref{sec:examples} presents the examples, with supporting calculations collected in the appendix. Section~\ref{sec:discussions} discusses some further directions of optimality, normalization and weighted likelihood envelopes, and connections with imprecise priors providing potential efficiency gains.

\section{Background}
\label{sec:background}

We briefly review the inferential model (IM) construction of \citet{martin2015inferential} in the possibility-contour formulation of \citet{martinhandbook}. Let $X\sim \mathrm P_\theta$, with
$\theta\in\Theta$, and suppose that the sampling model is represented
through the association
\[
    X=a(\theta,U), \qquad U\sim \mathrm P_U,
\]
where $u\in\mathcal U$  is a realization of auxiliary variable $U$, and the auxiliary distribution $\mathrm P_U$ is known. We assume that any relevant dimension reduction has already been carried
out, for example by working with a sufficient statistic or by conditioning
as in \citet[Section~3.2]{MartinConditional2015}. In the latter case,
the data determine a feature of the original auxiliary variable.
Conditioning on its observed value leaves only the remaining auxiliary
component to be predicted. We then use $X$ and $U$ for the reduced
variables and $\mathrm P_U$ for the conditional auxiliary law, with
sampling probabilities and expectations taken in the corresponding
conditional experiment. The smoothness, invertibility, and compatibility
conditions below apply to the reduced association. \citet{martinhandbook} illustrate auxiliary dimension reduction in their curved normal example. 

For fixed
$(x,\vartheta)$, write
\[
    \mathbb U_x(\vartheta)
    =
    \{u\in\mathbb U:x=a(\vartheta,u)\}.
\]
In the fully identified setting  considered below, this set is a
singleton, denoted $u_{x,\vartheta}$, and the association can also be
solved for $\vartheta=\theta_x(u)$. Throughout this smooth fully identified setting, after any required dimension reduction, assume $\mathcal U$ and $\Theta$ are $d$-dimensional Euclidean spaces (or suitable subspaces), and that the inverse maps are continuously differentiable with nonsingular Jacobians on the relevant interiors. We further assume that all auxiliary
values remain compatible with each observation, $\{u_{x,\vartheta}:\vartheta\in\Theta\}=\mathcal U$,
so that no part of the auxiliary space is ruled out by observing $x$.

\subsection{Predictive contours}
\label{sec:predictive-contours}

In the original IM formulation \citep{Martin2013}, uncertainty about the unobserved
auxiliary value is represented by a \emph{predictive random set}
$\mathcal S\subseteq\mathbb U$, with distribution
$\mathrm P_{\mathcal S}$. Combining $\mathcal S$ with the observation
$x$ gives the parameter-space random set
\[
    \Theta_x(\mathcal S)
    =
    \bigcup_{u\in\mathcal S}
    \{\vartheta:x=a(\vartheta,u)\}.
\]
Assuming this set is nonempty almost surely, the corresponding lower and
upper probabilities are
\[
    \lPi_x(A)
    =
    \mathrm P_{\mathcal S}
    \{\Theta_x(\mathcal S)\subseteq A\},
    \qquad
    \uPi_x(A)
    =
    \mathrm P_{\mathcal S}
    \{\Theta_x(\mathcal S)\cap A\neq\varnothing\},
\]
for assertions $A\subseteq\Theta$. They are dual, $\lPi_x(A)=1-\uPi_x(A^c)$.

The choice of predictive random set determines the efficiency of the
resulting inference. For a fixed assertion, admissibility arguments
motivate the use of nested predictive random sets, whose realizations are
ordered by inclusion \citep{MartinLiu2013Correction,martin2015inferential}.
For such a random set, define its hitting
function
\[
    \gamma(u)
    =
    \mathrm P_{\mathcal S}\{u\in\mathcal S\},
    \qquad u\in\mathbb U.
\]
Nesting implies that the auxiliary upper probability is completely
determined by $\gamma$ \citep{martinhandbook}:
\[
    \uPi^U(B)
    =
    \mathrm P_{\mathcal S}\{\mathcal S\cap B\neq\varnothing\}
    =
    \sup_{u\in B}\gamma(u).
\]
If $\sup_u\gamma(u)=1$, this is a possibility measure and $\gamma$ is
its \emph{possibility contour} \citep{DuboisPrade1988Possibility}. Its dual
$\lPi^U$ is a necessity measure, or equivalently a
consonant belief function.

Propagation through the association gives the parameter-space contour
\[
    \pi_x(\vartheta)
    =
    \sup_{u\in\mathbb U_x(\vartheta)}\gamma(u),
\]
which, in the fully identified case, reduces to $\pi_x(\vartheta)=\gamma(u_{x,\vartheta})$.
The complete inferential output is then obtained from
\[
    \uPi_x(A)
    =
    \sup_{\vartheta\in A}\pi_x(\vartheta),
    \qquad
    \lPi_x(A)
    =
    1-\uPi_x(A^c).
\]
Thus, once nesting is imposed, the predictive random set can be replaced
by the ordinary function $\gamma$ without changing the resulting IM.
Indeed, writing
\[
    S_\alpha
    =
    \{u:\gamma(u)\geq\alpha\},
    \qquad \alpha\in[0,1],
\]
the sets $S_\alpha$ are nested, and if
$V\sim\operatorname{Unif}(0,1)$ then
$\mathcal S=S_V$ has hitting function $\gamma$. We therefore work
directly with predictive contours in what follows.

The corresponding parameter-space cuts are $C_\alpha(x)=\{\vartheta:\pi_x(\vartheta)\geq\alpha\}$.
Under full identification, $C_\alpha(x)=\theta_x(S_\alpha)$, a relation that will be central to the later efficiency arguments.

\subsection{Auxiliary compatibility and sampling validity}\label{sec:aux-compat}

The predictive contour must be calibrated to the known auxiliary
distribution. This gives the following basic requirement.

\begin{definition}[Auxiliary compatibility]
\label{def:auxiliary-compatibility}
A contour $\gamma:\mathbb U\to[0,1]$ is \emph{compatible} with
$\mathrm P_U$ if
\[
    \mathrm P_U\{\gamma(U)<\alpha\}\leq\alpha,
    \qquad \alpha\in[0,1].
\]
More generally, an $x$-indexed family
$\{\gamma_x:x\in\mathbb X\}$ is auxiliary-compatible if
\[
    \mathrm P_U\{\gamma_x(U)<\alpha\}\leq\alpha
\]
for every fixed $x$ and every $\alpha\in[0,1]$.
\end{definition}

In terms of the contour cuts, compatibility is equivalently
\[
    \mathrm P_U(S_\alpha)\geq 1-\alpha,
    \qquad \alpha\in[0,1].
\]
We call the contour \emph{exact compatible} when
$\gamma(U)\sim\operatorname{Unif}(0,1)$, in which case
$\mathrm P_U(S_\alpha)=1-\alpha$ at every level.

Auxiliary compatibility should be distinguished from the
repeated-sampling calibration required of the induced parameter-space
contour.

\begin{definition}[Sampling validity]
\label{def:sampling-validity}
A family of parameter-space contours
$\pi=\{\pi_x:x\in\mathbb X\}$ is \emph{sampling valid} if
\[
    \sup_{\theta\in\Theta}
    \mathrm P_\theta\{\pi_X(\theta)<\alpha\}
    \leq\alpha,
    \qquad \forall \alpha\in[0,1].
\]
\end{definition}

With vacuous prior information, this is the contour form of \emph{strong validity} \citep{Martin2022a,Martin2023b}. For a single predictive contour, the next proposition gives the fully identified case of the predictive-random-set validity result of \citet{Martin2013}.

\begin{proposition}\label{prop:PRSvalidity}
Suppose $\pi_x(\vartheta)=\gamma(u_{x,\vartheta})$, where $\gamma$ is a single auxiliary contour compatible with
$\mathrm P_U$. If, under $\mathrm P_\theta$,
$u_{X,\theta}=U^\star$ for $U^\star\sim\mathrm P_U$, then $\pi$ is
sampling valid.
\end{proposition}

\begin{proof}
At the true parameter,
\[
    \pi_X(\theta)
    =
    \gamma(u_{X,\theta})
    =
    \gamma(U^\star).
\]
Hence
\[
    \mathrm P_\theta\{\pi_X(\theta)<\alpha\}
    =
    \mathrm P_U\{\gamma(U)<\alpha\}
    \leq\alpha,
\]
and the conclusion follows uniformly in $\theta$.
\end{proof}

Sampling validity transfers directly to assertion-wise inference. If
$\theta\in A$, then
$\uPi_X(A)\geq\pi_X(\theta)$, so
\[
    \sup_{\theta\in A}
    \mathrm P_\theta\{\uPi_X(A)<\alpha\}
    \leq\alpha.
\]
By duality, large lower probability assigned to a false assertion is
controlled in the corresponding way. In particular, the contour cuts
satisfy
\[
    \mathrm P_\theta\{\theta\in C_\alpha(X)\}
    \geq 1-\alpha.
\]

For a data-dependent family $\gamma_x$, pointwise auxiliary compatibility alone does not on its own imply sampling validity, because \[\pi_X(\theta)=\gamma_X(U^\star)=\gamma_{a(\theta, U^\star)}(U^\star),\] the same auxiliary realization can both select the contour and be evaluated by it. Sampling validity requires the stronger necessary diagonal condition: $P_U\{\gamma_{a(\theta, U)}(U)<\alpha\}\leq\alpha$ for all $\theta$ and all $\alpha$.

\subsection{Choice of predictive contour}
\label{sec:predictive-contour-choice}

Auxiliary compatibility does not determine a unique predictive contour.
A convenient general construction begins with a measurable score
$h:\mathbb U\to\mathbb R$, with larger values interpreted as more
plausible auxiliary values, and defines
\[
    \gamma_h(u)
    =
    \mathrm P_U\{h(U)\leq h(u)\}.
\]
If $h(U)$ has a continuous distribution, then the probability integral
transform gives
\[
    \gamma_h(U)\sim\operatorname{Unif}(0,1),
\]
so $\gamma_h$ is exact and, in particular, auxiliary-compatible. Thus
the score $h$ determines the ordering and geometry of the predictive
sets, while the probability transform supplies their calibration.

A distinguished choice is provided by the \emph{maximum-specificity}
principle from possibility theory \citep{DuboisEtAl2004}, applied to the auxiliary law by \citet{martinhandbook}. Among compatible contours that respect the density ordering, the principle selects the pointwise smallest one. Suppose that
$\mathrm P_U$ admits a density $f_U$ and that $f_U(U)$ has a continuous
distribution. The density-based maximum-specificity contour is
\begin{equation}\label{eq:max-spec}
    \gamma^{\mathrm{ms}}(u)
    =
    \mathrm P_U\{f_U(U)\leq f_U(u)\}.
\end{equation}
Hence auxiliary values are ordered according to their density under
$\mathrm P_U$, with high-density values receiving greater possibility.
Proposition~\ref{prop:PRSvalidity} implies sampling validity of the induced parameter-space
contour.

\paragraph{Example.} Let $Y_1,\ldots,Y_n$ be independent exponential
observations with mean $\theta>0$. The sufficient statistic
$X=\sum_{i=1}^n Y_i$ has association
\[
    X=\theta U, \qquad U\sim\operatorname{Gamma}(n,1),
\]
where the gamma distribution has shape $n$ and unit rate. The auxiliary
value compatible with $(x,\vartheta)$ is $u_{x,\vartheta}=x/\vartheta$.
Since $f_U(u)\propto u^{n-1}e^{-u}$, auxiliary maximum specificity
ranks candidates by $u^{n-1}e^{-u}$ evaluated at $u=x/\vartheta$.

\section{Optimal auxiliary-compatible contours}\label{sec:auxcontours} 

\citet[Section~4.5]{martin2015inferential} study theoretical optimality of predictive random sets, and \citet{Martin2025BvM} study asymptotic efficiency. \citet{raner2026efficiency} develops a decision-theoretic framework for finite-sample efficiency based on Choquet loss, the upper expectation of a penalty under the IM’s possibility measure. \emph{Contour concentration} corresponds to a particular choice of functional that penalizes the size of the contour's $\alpha$-cuts. This reduces to integrated size of cuts, or equivalently, the volume under the contour surface, and has close connections to classical notions of efficiency such as expected size of confidence sets.

The density-level sets of the maximum-specificity contour minimize auxiliary volume for a prescribed probability content \citep{NunezGarciaEtAl2003}. Contour concentration measures their volume after propagation to parameter space. The predictive set should thus ideally be chosen to make its image volume under the association small in parameter space, rather than merely small in auxiliary space. 

We will here impose contour-concentration risk functionals to guide the choice of efficient predictive auxiliary contour.

\subsection{An association-specific efficiency-ceiling}\label{sec:assoc-ceiling}

For fixed $x$, let $\gamma_x$ be an auxiliary-compatible contour satisfying the constraint \begin{equation}\label{eq:aux-comp}
    P_U\{\gamma_x(U)<\alpha\}\leq \alpha,
\end{equation} and let $S_\alpha(x)=\{u:\gamma_x(u)\geq \alpha \}$ be its auxiliary predictive $\alpha$-cut. The induced posterior $\alpha$-cut in parameter space is then, for the fully-identified case, \[C_\alpha(x)=\{\vartheta:u_{x,\vartheta}\in S_\alpha(x)\}=\theta_x(S_\alpha(x)).\] Let $\nu$ be the parameter-volume measure. In the default Euclidean setting, $\nu(d\vartheta)=d\vartheta$, the image-volume identity is \[\nu\{\theta_x(S_\alpha)\}=\int_{S_\alpha(x)}K_x(u)du,\qquad K_x(u)=|\operatorname{det}D_u\theta_x(u)|.\] More generally, if $\nu(d\vartheta)=v(\vartheta)d\vartheta$, replace $K_x(u)$ by $K_x^\nu(u)=v(\theta_x(u))|\operatorname{det}D_u\theta_x(u)|$. For notational simplicity, write $K_x(u)$ for the appropriate $\nu$-volume distortion.
Define the concentration loss 
\begin{equation}\label{eq:concloss}
    L(\gamma_x)=\int_0^1\nu\{C_\alpha(x)\}d\alpha.
\end{equation}

By Fubini/Tonelli, 
\begin{equation}\label{eq:layercake-aux}
    \begin{split}
        L_x(\gamma_x)&=\int_0^1\int_{\{\gamma_x(u)\geq\alpha\}}K_x(u)dud\alpha=\int_0^1\int_\mathbb U 1\{\gamma_x(u)\geq\alpha\}K_x(u) dud\alpha \\
        &=\int_\mathbb UK_x(u)\left[\int_0^11\{\gamma_x(u)\geq\alpha\}d\alpha\right] du=\int_\mathbb U K_x(u) \gamma_x(u)du.
    \end{split}
\end{equation}
Thus, for each fixed $x$, concentration minimization is an auxiliary-space weighted-volume problem. Moreover, as a consequence of nesting, both the loss \eqref{eq:concloss} and the compatibility constraint \eqref{eq:aux-comp} (equivalently $P_U\{S_\alpha(x)\}\geq1-\alpha$) decompose levelwise. At level $\alpha$, we must solve 
\begin{equation}
    \min_{S_\alpha\subseteq \mathbb U}\int_{S_\alpha(x)}K_x(u)du\qquad \text{subject to} \qquad \int_{S_\alpha(x)}f_U(u)du\geq1-\alpha.
\end{equation}
This is a Neyman--Pearson set-selection problem, with the auxiliary density $f_U$ playing the role of the probability to be retained and $K_x$ the volume cost. Indeed, define the parameter-volume measure $\mu_x(du)=K_x(u)du$, then the problem is \[\min_S \mu_x(S) \qquad \text{subject to} \qquad P_U(S)\geq 1-\alpha.\] The threshold form follows from the Neyman--Pearson lemma \citep[Theorem~3.2.1]{lehmann2005testing}. \citet{ScottNowak2006} give the corresponding minimum-volume formulation for a general reference measure. Consequently, the optimal sets are upper level sets of the density of $P_U$ relative to $\mu_x$: \[\frac{dP_U}{d\mu_x}(u)=\frac{f_U(u)}{K_x(u)}=:\rho_x(u).\]

\begin{theorem}[Fixed-data efficiency ceiling]
\label{thm:fix-data-ceiling}
Suppose the association is smooth and fully identified. For fixed $x$, let
$\rho_x(u)=f_U(u)/K_x(u)$, and assume that the distribution of $\rho_x(U)$,
$U\sim P_U$, is continuous. Define $\gamma_x^{\mathrm{ceil}}(u)
    =
    P_U\{\rho_x(U)\leq \rho_x(u)\}$, with induced parameter-space contour
\begin{equation}\label{eq:piceil}
    \pi_x^{\mathrm{ceil}}(\vartheta)
    =
    \gamma_x^{\mathrm{ceil}}(u_{x,\vartheta})
    =
    P_U\{\rho_x(U)\leq \rho_x(u_{x,\vartheta})\}.
\end{equation}
Then $\gamma_x^{\mathrm{ceil}}$ is auxiliary-compatible and, for every
$\alpha\in[0,1]$, its predictive cut 
\[S_\alpha^{\mathrm{ceil}}(x)=\{u:\gamma_x^{\mathrm{ceil}}(u)\geq\alpha\}\]
solves the levelwise weighted-volume problem above. Consequently,
$\gamma_x^{\mathrm{ceil}}$ minimizes the concentration loss
$L_x(\gamma_x)$ in (\ref{eq:concloss}) among all auxiliary-compatible
contours for the fixed data value $x$.
\end{theorem}

\begin{proof}
By \ref{eq:layercake-aux} and auxiliary compatibility
\ref{eq:aux-comp}, the concentration problem decomposes levelwise into
the Neyman--Pearson set-selection problem displayed immediately above the
theorem. Hence, at each $\alpha$, the minimizing set is an upper level set
of the density ratio $\rho_x=dP_U/d\mu_x=f_U/K_x$:
\[
    S_{\alpha,x}^\star
    =
    \{u:\rho_x(u)\geq c_{\alpha,x}\},
    \qquad
    P_U(S_{\alpha,x}^\star)=1-\alpha.
\]
By continuity of the distribution of $\rho_x(U)$, the probability integral
transform gives
$\gamma_x^{\mathrm{ceil}}(U)\sim\operatorname{Unif}(0,1)$, so
$\gamma_x^{\mathrm{ceil}}$ is auxiliary-compatible. Its $\alpha$-cuts are
precisely the sets $S_{\alpha,x}^\star$ up to null boundary sets. Thus it
minimizes the weighted volume at every level $\alpha$, and integrating
these levelwise inequalities in \ref{eq:layercake-aux} establishes
minimality of $L_x(\gamma_x)$.
\end{proof}

Note that the contour in \eqref{eq:piceil}, whilst auxiliary-compatible, \emph{is not in general sampling valid}. The ceiling contour uses an $x$-dependent auxiliary ranking. At the true value, $X=a(\theta, U^\star)$, but the score used to rank $U^\star$ is $\rho_X(\cdot)$, which depends on the same random $X$. Thus Proposition \ref{prop:PRSvalidity} does not apply, and validity is not guaranteed to hold in general. The contour acts as an efficiency ceiling and benchmark for any valid auxiliary-compatible contour.

\subsection{Averaged concentration risk}
The efficiency ceiling for auxiliary-compatible contours is determined by the fixed-$x$ Jacobian $K_x(u)$. But it is also this data-dependence that generally generally breaks the sampling validity of the resulting possibility measure in \eqref{eq:piceil}. To find an optimal valid procedure, we need to relax the optimization criteria. Interest is ultimately in concentration-risk. However, minimizing risk separately at a specified true parameter would generally produce a $\theta$-dependent, so attaining the pointwise optimum at every parameter value would require an oracle choice of auxiliary contour. To obtain a single implementable contour, we instead propose a risk-averaging over the parameter space via some design measure $w$ on $\Theta$. The measure $w$ is not a Bayesian prior; it is a risk-weighting device specifying how expected concentration loss is averaged over true values, creating a designed efficiency problem. The choice of $w$ specifies where precision is to be
prioritized. For example, a point mass targets risk at a prespecified
reference value, while a uniform design targets average risk over a
scientifically relevant range. In either case, validity holds throughout
$\Theta$.

Let $\gamma$ denote a single $P_U$ auxiliary-compatible contour.
The corresponding $w$-averaged concentration risk is
\begin{equation}\label{eq:w-risk}
    R_w(\gamma)=\int_\Theta E_\eta[L_X(\gamma)]w(d\eta).
\end{equation}
Using \eqref{eq:layercake-aux} and Tonelli, 
\begin{align*}
    R_w(\gamma)&=\int_\Theta E_\eta[L_X(\gamma)]w(d\eta)\\
    &=\int_\Theta E_\eta \left[\int_\mathbb U K_X(u)\gamma(u)du\right]w(d\eta)\\
    &=\int_\mathbb U \bar K_w(u)\gamma(u)du,
\end{align*}
where \begin{equation}\label{eq:Kbar}\bar K_w(u)=\int_\Theta E_\eta\left[K_X(u)\right]w(d\eta).\end{equation} As such we are just replacing the data-dependent Jacobian by its expected and $w$-averaged counterpart; a long-run volume price rather than observed volume price. The main issue is whether the integral in \eqref{eq:w-risk} exists or not.

For the derivations in this subsection, assume that the quantities above are well-defined for the contours under comparison. In particular, one may read the arguments first in the compact-$\Theta$ case, or more generally with a fixed finite design measure $w$ such that $\bar K_w(u)<\infty$ almost everywhere and the risks are finite. Appendix~\ref{app:noncompact-design} discusses two uses of compact
restriction on a non-compact parameter space. One may approximate
an infinite design measure by an increasing sequence of regions, or
choose a fixed, scientifically relevant region $C$ over which to
optimize average concentration risk. The latter choice specifies
where precision is to be prioritized, without assuming that
$\theta\in C$, and validity continues to hold throughout $\Theta$.

Rewriting the risk criterion in (\ref{eq:w-risk}) via layer-cake as before, the optimization problem is now 
\begin{equation}\label{eq:w-risk-layercake}\min_{\gamma}
    \int_0^1
    \int_{S_\alpha^\gamma}\bar K_w(u)\,du\,d\alpha\quad \text{subject to} \quad P_U(S_\alpha^\gamma)\ge 1-\alpha,
    \qquad
    \alpha\in[0,1].\end{equation}
The objective and constraint still decompose by level, and the levelwise problem is
\begin{equation}\label{eq:w-levelwise-problem}
    \min_{S_\alpha\subseteq\mathbb U}
    \int_{S_\alpha}\bar K_w(u)\,du
    \quad
    \text{subject to}
    \quad
    \int_{S_\alpha}f_U(u)\,du\ge 1-\alpha.
\end{equation}
Consequently, by the same Neyman--Pearson argument, the optimal cuts are given by thresholding of \[\rho_w(u)=\frac{f_U(u)}{\bar K_w(u)}.\]

\begin{theorem}[Optimal $w$-averaged single PRS]
\label{thm:prs-w-choquet-ceiling}
Suppose the association is smooth and fully identified, and let $w$ be a
design measure such that the $w$-averaged risk in \ref{eq:w-risk} is finite
for the contours under consideration. Define
\[
\bar K_w(u)
=
\int_\Theta E_\eta\{K_X(u)\}\,w(d\eta),
\qquad
\rho_w(u)
=
\frac{f_U(u)}{\bar K_w(u)},
\]
and assume that $0<\bar K_w(u)<\infty$ for $P_U$-almost every $u$, and that the distribution of $\rho_w(U)$,
$U\sim P_U$, is continuous. Then the single auxiliary contour
\begin{equation}\label{eq:prs-w-optimal-aux-contour}
\gamma_w^\star(u)
=
P_U\{\rho_w(U)\leq\rho_w(u)\}
\end{equation}
minimizes the $w$-averaged concentration risk $R_w$ in \ref{eq:w-risk}
among all single auxiliary contours compatible with $P_U$.

Equivalently, for every $\alpha\in[0,1]$, its predictive cut is
$S_{\alpha,w}^\star
=
\{u:\rho_w(u)\geq c_{\alpha,w}\}$, where $c_{\alpha,w}$ is chosen so that $P_U(S_{\alpha,w}^\star)=1-\alpha$.
The induced parameter-space contour
\[
\pi_{x,w}^\star(\vartheta)
=
\gamma_w^\star(u_{x,\vartheta})
\]
is sampling valid. If, moreover, $\{u_{x,\vartheta}:\vartheta\in\Theta\}
=
\mathbb U$
for every $x$, then it is normalized:
\[
\sup_{\vartheta\in\Theta}
\pi_{x,w}^\star(\vartheta)
=
1.
\]
\end{theorem}

\begin{proof}
By the risk representation in \ref{eq:w-risk-layercake}, the objective and
compatibility constraint decompose levelwise, with each level governed by
\ref{eq:w-levelwise-problem}. This is exactly the fixed-data problem of
Theorem \ref{thm:fix-data-ceiling}, with $K_x$ replaced by $\bar K_w$.
Hence the same Neyman--Pearson argument gives the minimizing cuts
\[
    S_{\alpha,w}^{\star}
    =
    \{u:\rho_w(u)\geq c_{\alpha,w}\},
    \qquad
    \rho_w(u)=\frac{f_U(u)}{\bar K_w(u)},
\]
where $c_{\alpha,w}$ is chosen so that $P_U(S_{\alpha,w}^{\star})=1-\alpha$. Since these are upper level sets of the common score $\rho_w$, they are
nested and generate the contour in \ref{eq:prs-w-optimal-aux-contour}.
Levelwise optimality therefore implies optimality of the integrated
$w$-averaged risk.

By continuity of the distribution of $\rho_w(U)$, the probability integral
transform gives 
\[
    \gamma_w^\star(U)\sim\operatorname{Unif}(0,1),
\]
so $\gamma_w^\star$ is compatible with $P_U$. Sampling validity of the
induced contour $\pi^{PRS,w}$ then follows directly from Proposition
\ref{prop:PRSvalidity}.

Finally, in the regular fully-identified case, the inverse association
ranges over $\mathbb U$ for every fixed $x$. Therefore
\[
    \sup_{\vartheta\in\Theta}\pi_x^{PRS,w}(\vartheta)
    =
    \sup_{u\in\mathbb U}\gamma_w^\star(u)
    =
    1,
\]
which proves normalization.
\end{proof}

\subsection{Observed and design-averaged volume costs}
The fixed-data ceiling contour in Theorem~\ref{thm:fix-data-ceiling} and the design-averaged optimal contour in Theorem~\ref{thm:prs-w-choquet-ceiling} arise from the same levelwise set-selection problem, but with different notions of volume cost. For fixed $x$, the relevant cost of including $u$ in a predictive set is the \emph{observed} pullback volume $K_x(u)$ giving the ordering $\rho_x(u)=f_U(u)/K_x(u)$. By contrast, requiring a single auxiliary contour to be chosen independently of the realized data replaces this observed cost by its design-averaged counterpart $\bar K_w(u)$ in (\ref{eq:Kbar}), and hence gives the ordering $\rho_w(u)=f_U(u)/\bar K_w(u)$.

Thus, $K_x(u)$ and $\bar K_w(u)$ can be viewed as two different volume prices. The former measures the parameter-space volume induced by the auxiliary point $u$ for the actually observed $x$, whereas the latter measures its expected volume cost over repeated sampling and over the design measure $w$. We have,
\begin{equation}\label{eq:design-score}
    \rho_w(u)=\rho_x(u)\frac{K_x(u)}{\bar K_w(u)}.
\end{equation}
The factor $K_x(u)/\bar K_w(u)$ therefore measures the discrepancy between the observed and long-run volume costs assigned to $u$. At the fixed-data level there is no penalty for adapting the auxiliary ordering to the realized pullback geometry, while a single auxiliary contour must use one ordering for all possible samples. Going from $K_x$ to $\bar K_w$ is precisely the cost of imposing this data-independent auxiliary geometry.

There are, however, important cases in which this distinction disappears. Suppose for example that the pullback-volume factor separates  $K_x(u)=c(x)K_0(u)$, for some positive functions $c$ and $K_0$. Then \[\rho_x(u)=\frac{1}{c(x)}\frac{f_U(u)}{K_0(u)},\] so the ordering of the auxiliary space induced by $\rho_x$ is independent of $x$. Moreover \[\bar K_w(u)=K_0(u)\int_\Theta E_\eta[c(X)]w(d\eta), \] whenever the final integral is finite, and consequently \[\rho_w(u) \propto \frac{f_U(u)}{K_0(u)}.\] Hence the fixed-data and design-averaged constructions induce the same nested family of predictive sets.
\begin{corollary}\label{cor:factoriz}
    Suppose the assumptions of Theorem~\ref{thm:fix-data-ceiling} and Theorem~\ref{thm:prs-w-choquet-ceiling} hold, and that the observed parameter-volume distortion factors as \[K_x(u)=c(x)K_0(u),\]
    for positive functions $c$ and $K_0$. Then the fixed-data ceiling auxiliary contour is independent of $x$ and coincides with the design-averaged optimal contour: $\gamma_x^{\mathrm{ceil}}=\gamma_w^\star$. This common contour is the same for every design measure $w$ satisfying the assumptions of Theorem~\ref{thm:prs-w-choquet-ceiling}. Consequently, the induced ceiling contour is sampling valid.
\end{corollary}
\begin{proof}
    By the factorization and the assumption on $w$, \[\bar K_w(u)=A_wK_0(u), \qquad A_w:=\int_\Theta E_\eta[c(X)]w(d\eta) \in (0,\infty).\] Therefore, \[\rho_x(u)=\frac{1}{c(x)}\frac{f_u(u)}{K_0(u)}, \qquad \rho_w(u)=\frac{1}{A_w}\frac{f_U(u)}{K_0(u)}.\] Positive rescaling leaves the auxiliary ordering and its probability transform unchanged, so \[\gamma_x^{\mathrm{ceil}}(u)=P_U\left\{\frac{f_U(U)}{K_0(U)}\le \frac{f_U(u)}{K_0(u)} \right\} = \gamma_w^\star (u).\] Sampling validity now follows from Theorem~\ref{thm:prs-w-choquet-ceiling}.
\end{proof}
In particular, under this factorization the fixed-data ceiling is no longer merely an oracle benchmark. Its auxiliary ordering can itself be represented by a single auxiliary contour and therefore inherits the corresponding sampling-validity property.

\paragraph{Example (continued).} In the exponential example, $\theta_x(u)=x/u$ and $E_\eta(X)=n\eta$,
so the observed and design-averaged costs for ordinary parameter
length are
\[
    K_x(u)=\frac{x}{u^2},
    \qquad
    \bar K_w(u)=\frac{n}{u^2}\int_0^\infty \eta\,w(d\eta).
\]
For any admissible $w$, the two costs differ only by a factor
independent of $u$. Thus Corollary~\ref{cor:factoriz} gives the common
ordering $u^2 f_U(u)\propto u^{n+1}e^{-u}$. The extra factor $u^2$
favors larger auxiliary values, where a given auxiliary length maps
to a shorter parameter interval. Since the ordering is independent
of $x$, the ceiling itself is sampling valid.

This factorization is a convenient sufficient condition rather than a necessity for the phenomenon. What matters is that the data dependence of the fixed-$x$ score does not alter the ordering of the auxiliary space. More generally, whenever the family $\{\rho_x:x\in \mathrm X\}$ has a common auxiliary ordering, the apparent data dependence of the fixed-data optimization is irrelevant for the auxiliary ranking. With ordinary parameter length, additive location $x=\theta+u$ gives $K_x(u)=1$, and positive scale $x=\theta u$ gives $K_x(u)=x/u^2$. Both therefore satisfy the factorization.

\section{Sampling-valid repairs of the efficient ceiling contour}\label{sec:samplingrepair}
The efficient ceiling contour in Theorem~\ref{thm:fix-data-ceiling} uses the observed-data optimal ranking \[\rho_x(u)=\frac{f_U(u)}{K_x(u)},\] and minimizes volume for the observed $x$. Because this ranking is data-dependent, the resulting $x$-indexed auxiliary contour is generally not sampling valid. The construction of the design-averaged contour in Theorem~\ref{thm:prs-w-choquet-ceiling} is essentially motivated by rescuing this construction with a single compatible predictive auxiliary contour from which validity follows. This is achieved by replacing the volume cost $K_x$ by the design-averaged volume cost $K_w$, which effectively removes the data-dependence from the score, but then changes the efficiency problem into one of long-run integrated volume minimization instead.

Proposition~\ref{prop:PRSvalidity} is just one strategy to establish validity however. We do not require a valid contour to be of form $\pi_x(\vartheta)=\gamma(u_{x,\vartheta})$ for one common $\gamma$. The aim is not compatibility, but validity in parameter space: \[\sup_\theta P_\theta\{\pi_X(\theta)\leq\alpha\}\leq\alpha, \qquad \alpha\in[0,1].\] Allowing for candidate-specific calibration: $\pi_x(\vartheta)=\gamma_\vartheta(u_{x,\vartheta})$ — might allow for more flexible means of validity repair. 

In this section the focus is now on the goal: retain as much as possible the fixed-data ceiling's local geometry; repair the failure of sampling validity without forcing a single compatible auxiliary ordering across all datasets. There are two immediate repair strategies: repair the score $\rho_x(u)$ underlying the ceiling, or repair the ceiling contour $\pi_x$ after it has already been constructed. Crucial to this end is a transform that allows for candidatewise validification. \citet[Section~2.1]{Martin2015Plausibility} uses the sampling distribution of a candidate-specific score for this purpose.

\subsection{Candidatewise validification}
For any candidate-dependent score $T_x(\vartheta)$, its sampling distribution is given by $F_\vartheta^T(t)=P_\vartheta\{T_X(\vartheta)\leq t\}$. We define the candidate-wise validification of the score $T_x(\vartheta)$ by 
\begin{equation}\label{eq:candidate-validification}
    \mathcal V[T]_x(\vartheta):=F_\vartheta^T\{T_x(\vartheta)\}=P_\vartheta\{T_X(\vartheta)\leq T_x(\vartheta)\}.
\end{equation}
Then, by properties of the probability-integral-transform (PIT), we have \[P_\theta\{\mathcal V[T]_X(\theta)\leq\alpha\}\leq\alpha,\] and under continuity, $\mathcal V[T]_X(\theta)\sim \operatorname{Unif}(0,1)$.

Moreover, this PIT validification is the pointwise smallest sampling-valid non-decreasing transformation of the chosen scalar score. That is, for any valid $\pi_x(\vartheta)=g_\vartheta\{T_x(\vartheta)\}$, with $g_\vartheta$ non-decreasing, validity requires $g_\vartheta(t)\geq F_\vartheta^T(t)$. Indeed, with $T=T_X(\vartheta)$, monotonicity and sampling validity give
\[
  F_\vartheta^T(t)
    \leq P_\vartheta\{g_\vartheta(T)\leq g_\vartheta(t)\}
    \leq g_\vartheta(t).
\]

\subsection{Localized auxiliary repair}\label{sec:loc-repair}
For each candidate value $\vartheta$, consider varying the auxiliary value $u$. This traces out the data that could have been observed if $\vartheta$ were true:  $x_\vartheta(u)=a(\vartheta, u).$ Since the association is assumed fully identified, this is well defined. The fixed-data efficiency score in Theorem~\ref{thm:fix-data-ceiling} can then be evaluated at each of those hypothetical datasets. At each $x_\vartheta(u)$, we ask how volume would be distorted by that hypothetical observation: we define the \emph{local image-volume Jacobian} $K_\vartheta^{\mathrm {loc}}(u):=K_{x_\vartheta(u)}(u)$,
and the corresponding local score 
\begin{equation}
    \rho_\vartheta^{\mathrm{loc}}(u)=\frac{f_U(u)}{K_\vartheta^{\mathrm{loc}}(u)}.
\end{equation}
Rather than fixing $x$ and using the volume-optimal ranking for that single data value, localization lets the geometry evolve along the datasets that $\vartheta$ could generate.

Transforming the score as $\gamma_\vartheta(u)=P_U\{\rho_\vartheta^{\mathrm{loc}}(U)\leq \rho_\vartheta^{\mathrm{loc}}(u)\}$ yields a valid p-value function, but it need not yield a normalized possibility contour, even though for fixed $\vartheta$ this local auxiliary contour is normalized on $\mathrm U$. Because the contour's maximizer itself moves with $\vartheta$, the maximizer is not guaranteed to be compatible with a pair $(x,\vartheta)$ for the observed $x$.

Once $x$ is actually observed, the pair $(x, \vartheta)$ identifies the unique compatible auxiliary value $u_{x,\vartheta}$. Crucially, \[x_\vartheta(u_{x,\vartheta})=a(\vartheta,u_{x,\vartheta})=x.\] Therefore, at the particular auxiliary point compatible with the observation, the localized score reduces back to the observed-$x$ efficiency score: \[\rho_\vartheta^{\mathrm{loc}}(u_{x,\vartheta})=\rho_{a(\vartheta, u_{x,\vartheta})}(u_{x,\vartheta})=\rho_x(u_{x,\vartheta}).\]
For every $x\in\mathcal X$, assume that \[0<M(x):=\sup_{v\in\mathcal U}\rho_x(v)<\infty\] and that the supremum is attained at some $\hat v_x\in \mathcal U$.
Define the normalized score
\begin{equation} \tilde\rho_\vartheta(u)=\frac{\rho_{a(\vartheta, u)}(u)}{M(a(\vartheta,u))}.
\end{equation}
Now each hypothetical observation $x_\vartheta(u)$ is assessed relative to its own fixed-data efficiency ceiling. At the actually observed pair $(x, \vartheta)$, we have $\tilde\rho_\vartheta(u_{x,\vartheta})=\rho_x(u_{x,\vartheta})/M(x)$. The auxiliary contour given by this normalized score is then $\tilde\gamma_\vartheta(u)=P_U\{\tilde\rho_\vartheta(U)\leq \tilde\rho_\vartheta(u)\}$, thus the final localized contour is 
\begin{equation}\label{eq:pi-loc}
    \pi_x^{\mathrm{loc}}(\vartheta)=\tilde \gamma_\vartheta(u_{x,\vartheta}).
\end{equation}

Let $\tau_x(\vartheta)=\rho_x(u_{x,\vartheta})/M(x)$. The localized contour in (\ref{eq:pi-loc}) is then equivalently given as \[P_U\left\{\frac{\rho_{a(\vartheta, U)}(U)}{M(a(\vartheta,U))}\leq \frac{\rho_x(u_{x,\vartheta})}{M(x)}\right\}=P_\vartheta\{\tau_X(\vartheta)\leq\tau_x(\vartheta)\},\] and by (\ref{eq:candidate-validification}) this is valid. In particular, set $\hat \vartheta_x=\theta_x(\hat v_x)$. Then $\tau_x(\hat \vartheta_x)=1$. Since $\tau_X(\hat \vartheta_x)\le 1$,  \[\pi_x^{\mathrm{loc}}(\hat \vartheta_x)=P_{\hat \vartheta_x}\{\tau_X(\hat \vartheta_x)\leq 1\}=1.\] Thus the validified localized contour is normalized.

The candidate-indexed auxiliary contours $\tilde\gamma_\vartheta$ are combined by diagonal evaluation. At $\vartheta$ we use $\tilde\gamma_\vartheta(u_{x,\vartheta})$, and then take the supremum over the assertion. This is the nested case of the fusion construction of \citet[Section~4.2]{Martin2021}. The attained maximum above verifies its normalization condition.

\subsubsection{Sample-space localization and the LR IM}\label{sec:validified-LR}
There is a closely related localization which leads to the familiar likelihood-ratio IM. Recall that the maximum-specificity construction on $\mathrm U$ ranks auxiliary values by $f_U(u)$, corresponding to minimum auxiliary-space volume of contour cuts for a prescribed $P_U$-probability. As before, rather than retaining auxiliary volume itself, fix a candidate $\vartheta$ but now consider the image of a predictive set under $x_\vartheta(u)=a(\vartheta,u)$ in sample-space. Let \begin{equation}
    \label{eq:samp-change-volume}
    A_\vartheta(u)=|\operatorname{det}D_ua(\vartheta,u)|
\end{equation} denote the corresponding sample-space image-volume factor. Then, for $S\subseteq \mathrm U $, the candidate-local analogue of the maximum-specificity problem is therefore \[\min_{S\subseteq \mathrm U}\int_SA_\vartheta(u)du \qquad \text{subject to} \qquad P_U(S)\geq1-\alpha.\] 

By the same Neyman--Pearson argument as before, the optimal sets are upper level sets of $\rho_\vartheta^X(u)=f_U(u)/A_\vartheta(u)$. But a change of variables through the association gives \begin{equation}\label{eq:samp-changeof-var}
    \frac{f_U(u)}{A_\vartheta(u)}=p_\vartheta\{a(\vartheta,u)\}.
\end{equation} Thus the likelihood emerges as the score obtained by localizing the auxiliary maximum-specificity construction in the $U \rightarrow \mathcal X$ direction. In particular, at the auxiliary value $u_{x,\vartheta}$ compatible with the observed pair $(x, \vartheta)$, we have $\rho_\vartheta^X(u_{x, \vartheta})=p_\vartheta(x)$.

As in Section~\ref{sec:loc-repair}, direct candidatewise validification of this local score need not produce a normalized contour over $\vartheta$. We therefore normalize each hypothetical observation by its maximal score. Assume likewise that, for every $x$, the likelihood supremum $\sup_\eta p_\eta (x)$ is finite, strictly positive, and attained. Define \[\tilde \rho _\vartheta^X(u)=\frac{p_\vartheta\{a(\vartheta,u)\}}{\sup_{\eta\in\Theta}p_\eta\{a(\vartheta,u)\}}.\] At the observed compatible point this reduces to the familiar relative likelihood \[R_x(\vartheta)=\tilde \rho _\vartheta^X(u_{x,\vartheta})=\frac{p_\vartheta(x)}{\sup_{\eta\in\Theta}p_\eta\{x\}}.\] Candidatewise validification now gives 
\begin{equation}
    \pi_x^{\mathrm{LR}}(\vartheta)=P_U\{\tilde \rho _\vartheta^X(U)\le \tilde \rho _\vartheta^X(u_{x,\vartheta})\},
\end{equation}
which is the usual likelihood-ratio IM \citep{Martin2015Plausibility}. The normalization changes the original minimum-Lebesgue-volume problem; the relative likelihood instead induces a weighted sample-space ordering, but preserves the candidate-local interpretation.

This construction provides a useful contrast with the localization in Section~\ref{sec:loc-repair}. Both begin from an auxiliary Neyman--Pearson ordering, then localize, normalize, and validify candidatewise. The difference is the direction in which auxiliary volume is transported. The likelihood-ratio construction localizes through $u\mapsto x_\vartheta(u)$ and is therefore naturally targeted towards specificity in sample-space. Section~\ref{sec:loc-repair} instead begins with the fixed-data ordering induced by $u\mapsto \theta_x(u)$, whose local cost is parameter-space volume. Since the concentration loss in (\ref{eq:concloss}) is itself defined in terms of parameter-space volume, the latter localization is more directly aligned with the efficiency criterion. This does not imply a general dominance after normalization and validification, but it provides a structural reason to expect the localized parameter-volume construction to tend towards smaller concentration loss than the ordinary LR IM.

\subsection{Direct repair of the efficiency ceiling}\label{sec:direct-rep}
Here we propose a more direct method of repair. Recall the fixed-data optimal contour from Theorem~\ref{thm:fix-data-ceiling}: $\pi_x^{\mathrm{ceil}}(\vartheta)=P_U\{\rho_x(U)\le \rho_x(u_{x,\vartheta})\}$. Instead of repairing the score $\rho_x$, we can take the contour $\pi_x^{\mathrm{ceil}}$ itself as the object whose order should be preserved as much as possible, and calibrate that contour directly in sampling space. That is, we validify the ceiling contour as
\begin{equation}\label{eq:validified-ceiling}
    \pi_x^{Vc}(\vartheta)=\mathcal{V}[\pi^{\mathrm{ceil}}]_x(\vartheta)=P_\vartheta\{\pi_X^{\mathrm{ceil}}(\vartheta)\le \pi_x^{\mathrm{ceil}}(\vartheta)\}.
\end{equation}

Since $\pi_x^{\mathrm{ceil}}(\vartheta)$ is already normalized before validification, $\pi_x^{\mathrm{Vc}}(\vartheta)$ is normalized by the same argument as in Section~\ref{sec:loc-repair}. Validity follows from (\ref{eq:candidate-validification}).

\section{Fiducial Connections}\label{sec:connections}
The constructions in Sections \ref{sec:auxcontours} and \ref{sec:samplingrepair} also admit a common fiducial interpretation. In the smooth fully identified setting, the inverse association determines a generalized fiducial distribution whose density generates the efficiency ceiling. The sampling-valid constructions correspond to different treatments of this density ordering and its calibration. 

\subsection{The generalized fiducial distribution and the efficient ceiling}
Assume that the inverse maps are one-to-one and continuously differentiable on their supports, with all auxiliary values compatible. For fixed $x$, denote by $\mathrm R_x$ the distribution obtained by propagating $P_U$ through the inverse association: \[\Theta_x=\theta_x(U'),\qquad U'\sim P_U, \qquad \mathrm R_x(B)=P_U\{u:\theta_x(u)\in B\}.\] In this setting, $\mathrm R_x$ is the generalized fiducial distribution \citep[e.g.,][]{Hannig2016}. 

Retain the parameter-volume measure $\nu(d\vartheta)=v(\vartheta)d\vartheta$ and the convention from Section~\ref{sec:assoc-ceiling} that $K_x$ includes its density $v$. The likelihood-times-Jacobian expression in \eqref{eq:fid-density} specializes the general fiducial density formula of \citet[Theorem~1]{Hannig2016} to the invertible case and expresses it relative to $\nu$.
\begin{proposition}\label{prop:fid-ceiling}
    Suppose the parameter and auxiliary derivatives of the association are nonsingular. Let $J(x,\vartheta)=|\operatorname{det}D_\vartheta a(\vartheta,u_{x,\vartheta})|$. Then 
    \begin{equation}
        \label{eq:fid-density}
        r_x^\nu(\vartheta):=\frac{d\mathrm R_x}{d\nu}(\vartheta)=\rho_x(u_{x,\vartheta})=\frac{p_\vartheta(x)J(x,\vartheta)}{v(\vartheta)}.
    \end{equation} Consequently, the efficient ceiling is the density-based maximum-specificity probability-to-possibility transform \begin{equation}\label{eq:fid-ceiling}
        \pi_x^{\mathrm{ceil}}(\vartheta)=\mathrm R_x\{r_x^\nu(\Theta)\leq r_x^\nu(\vartheta)\}.
    \end{equation}
\end{proposition}
\begin{proof}
    Since $\theta_x$ and $u_x:\vartheta \mapsto u_{x,\vartheta}$ are inverses, change of variables gives \[r_x^\nu(\vartheta)=\frac{f_U(u_{x,\vartheta})}{v(\vartheta)|\operatorname{det}D_u\theta_x(u_{x,\vartheta})|}=\frac{f_U(u_{x,\vartheta})}{K_x(u_{x,\vartheta})}.\] Implicit differentiation of $a(\theta_x(u),u)=x$ with respect to $u$ gives \[D_u\theta_x(u)=-\bigl[D_\vartheta a(\theta_x(u),u)\bigr]^{-1}D_ua(\theta_x(u),u).\] Thus $K_x(u_{x,\vartheta})=v(\vartheta)A_\vartheta(u_{x,\vartheta})/J(x,\vartheta)$, with $A_\vartheta(u)$ given as in (\ref{eq:samp-change-volume}). Moreover, change of variables gives (\ref{eq:samp-changeof-var}), which establishes the density identity. Finally, since $u_x$ and $\theta_x$ are inverse maps, we have $u_{x,\Theta_x}=U'$ and $r_x^\nu(\Theta_x)=\rho_x(U')$. Consequently, substitution into (\ref{eq:piceil}) gives the representation in~(\ref{eq:fid-ceiling}).
\end{proof}

For Lebesgue parameter volume and a continuous density rank, \citet{martinhandbook} give \eqref{eq:fid-ceiling} by applying the probability-to-possibility transform to the fiducial density. For the general form of this transform, see \citet{reIM}. Theorem~\ref{thm:fix-data-ceiling} identifies its minimum-volume property and extends the construction to a specified measure $\nu$.

Write $s_x(\vartheta)=r_x^\nu(\vartheta)$, and call this the \emph{association-induced fiducial likelihood}, relative to $\nu$. Proposition~\ref{prop:fid-ceiling} identifies the fiducial density with the parameter-side score selected by the concentration problem. Under the continuity assumption of Theorem~\ref{thm:fix-data-ceiling}, the ceiling cuts are therefore highest-density regions of $\mathrm R_x$, with minimum $\nu$-volume at their prescribed fiducial content. Moreover, $\mathrm R_x$ is a \emph{saturating inner probability} of the possibility measure $\Pi_x^{\mathrm{ceil}}(B)=\sup_{\vartheta \in B}\pi_x^{\mathrm{ceil}}(\vartheta)$: it is bounded from above by $\Pi_x^{\mathrm{ceil}}(B)$ for every assertion $B$, and satisfies $\mathrm R_x\{C_\alpha^{\mathrm{ceil}}(x)\}=1-\alpha.$ 

This calibration is fixed-data, not repeated-sampling (for more on saturating inner measures of \emph{valid} IMs, see \citet{reIM}). After fixing $x$, fiducial inversion continues to use $P_U$ and propagates a fresh $U'$ through $\theta_x$, giving $\pi_x^{\mathrm{ceil}}(\Theta_x)\sim\operatorname{Unif}(0,1)$. This fixed-data calibration clarifies the sampling limitation of the fiducial \emph{continue-to-regard} argument \citep[p.~109]{Dempster1963}. The resulting uniformity concerns $\pi_x^{\mathrm{ceil}}\{\theta_x(U')\}=\gamma_x^{\mathrm{ceil}}(U')$. Under repeated sampling at the true parameter, however, we have \[X=a(\theta, U^\star), \qquad \pi_X^{\mathrm{ceil}}(\theta)=\gamma_{a(\theta, U^\star)}^{\mathrm{ceil}}(U^\star).\] Here the auxiliary realization being ranked also determines the ranking. Continuing to use $P_U$ after fixing $x$ defines the fiducial law, but does not establish the diagonal sampling-calibration condition of Section~\ref{sec:aux-compat}.

\subsection{Fiducial interpretations of the sampling-valid constructions}
The design-averaged construction transforms the same fiducial law using a different score. By (\ref{eq:design-score}), define \[s_{w,x}(\vartheta)=s_x(\vartheta)\frac{K_x(u_{x,\vartheta})}{\bar{K}_w(u_{x,\vartheta})}=\rho_w(u_{x,\vartheta}).\] Since under $\Theta\sim\mathrm R_x$, $s_{w,x}(\Theta)$ has the same distribution as $\rho_w(U)$, we get \[\pi^\star_{x,w}(\vartheta)=\mathrm R_x\{s_{w,x}(\Theta)\leq s_{w,x}(\vartheta)\}.\] Thus $\mathrm R_x$ also saturates the $w$-optimal contour. The design measure changes the ordering through the observed-to-design-averaged volume correction, rather than changing the underlying probability law.

The localized construction in Section~\ref{sec:loc-repair} has a direct fiducial representation. Assume that the suprema used below are finite and attained. Define the \emph{relative fiducial likelihood} 
\begin{equation}
    \label{eq:fid-likelihood}
    RF_x(\vartheta)=\frac{s_x(\vartheta)}{\sup_{\eta \in \Theta}s_x(\eta)}.
\end{equation}
Since $\sup_\eta s_x(\eta)=M(x)$, this is exactly the normalized score $\tau_x(\vartheta)$ in Section~\ref{sec:loc-repair}. Consequently, 
\begin{equation}
    \label{eq:validified-fid-like}
    \pi_x^{\mathrm{loc}}(\vartheta)=P_\vartheta\{RF_X(\vartheta)\le RF_x(\vartheta)\}.
\end{equation}
Thus localization amounts to candidatewise sampling validification of the relative fiducial likelihood.

By (\ref{eq:fid-density}), \[RF_x(\vartheta)=\frac{p_\vartheta(x)J(x,\vartheta)/v(\vartheta)}{\sup_{\eta\in\Theta} p_\eta(x)J(x,\eta)/v(\eta)}.\] The difference from the relative likelihood $R_x(\vartheta)$ in Section~\ref{sec:validified-LR} is therefore the association-induced correction $J(x,\vartheta)/v(\vartheta)$. In particular, if this correction is independent of $\vartheta$ for every $x$, then $RF_x=R_x$, and $\pi_x^{\mathrm{loc}}=\pi_x^{\mathrm{LR}}$.

\paragraph{Example (continued).} For the exponential association with ordinary parameter length, this
correction is $J(x,\vartheta)=x/\vartheta$. Hence
$r_x(\vartheta)=p_\vartheta(x)\,x/\vartheta$.
Writing $u=x/\vartheta$, the correction takes the ordinary likelihood
ordering $u^n e^{-u}$ to the fiducial ordering $u^{n+1}e^{-u}$
obtained from the parameter-volume criterion.

To separate fiducial calibration from sampling calibration, define \[H_x(t)=\mathrm R_x\{RF_x(\Theta) \le t\},\qquad F_\vartheta(t)=P_\vartheta\{RF_X(\vartheta)\le t\}.\] Then \[\pi_x^{\mathrm{ceil}}(\vartheta)=H_x\{RF_x(\vartheta)\},\qquad \pi_x^{\mathrm{loc}}(\vartheta)=F_\vartheta\{RF_x(\vartheta)\}.\] The direct repair in Section~\ref{sec:direct-rep} instead becomes \[\pi_x^{\mathrm{Vc}}(\vartheta)=P_\vartheta[H_X\{RF_X(\vartheta)\}\le H_x\{RF_x(\vartheta)\}].\] The two repairs therefore impose sampling calibration at different stages. Localization validifies the relative fiducial likelihood; direct repair first converts that score into its fixed-data fiducial rank, and then validifies the rank. 

\subsection{Further literature connections}
$\mathrm R_x$ being saturating is not specific to the optimizing scores. Inverse consistency of $u_x$ and $\theta_x$ gives \[\mathrm R_x=(\theta_x)_\#P_U,\qquad (u_x)_\#\mathrm R_x =P_U.\] Denote by $\Gamma_U=\{\gamma: \mathrm U \rightarrow [0,1]: \gamma(U)\sim \operatorname{Unif}(0,1)\}$ the class of exact auxiliary contours. Let $\mathcal R=\{R_x:x\in \mathrm X\}$ be the fiducial family induced by the inverse association. Its \emph{association-generated orbit} is \[\mathcal O_a(\mathcal R)=\{\pi^\gamma:\pi_x^\gamma(\vartheta)=\gamma(u_{x,\vartheta}), \quad\gamma \in \Gamma_U\}.\] Every member is exactly sampling-calibrated and saturated by the same fiducial law. That is, every cut $C_\alpha^\gamma(x)=\{\vartheta:\pi_x^\gamma(\vartheta)\ge \alpha\}$ satisfies \[\mathrm R_x\{C_\alpha^\gamma(x)\}=P_\theta\{\theta\in C_\alpha^\gamma(X)\}=1-\alpha.\] Different auxiliary orderings therefore produce different calibrated confidence set geometries with the same fiducial probability content. Theorem~\ref{thm:prs-w-choquet-ceiling} selects the specific efficient member $\gamma^\star_w$. The efficient ceiling is likewise saturated by $\mathrm R_x$, but generally lies outside this orbit because its auxiliary ordering depends on the observation. 

The orbit supplies many confidence-region families for which fiducial probability and sampling probability agree exactly. This suggests that a sufficiently stable perturbation of an orbital confidence region can retain both calibrations approximately. It suffices that the symmetric difference between the two regions has small fiducial probability at the observed $x$ and rarely contains the true $\theta$ under repeated sampling. This may provide a partial structural explanation for the often good approximate accuracy of fiducial-derived confidence regions. For examples and an asymptotic coverage result, see \citet{Hannig2016}.

Connections between fiducial inference and IMs have been established from several directions. \citet{martinhandbook} interpret the primitive fiducial distribution as an efficiency ceiling for predictive random set IMs, drawing on the bound on belief established by \citet[Proposition~1]{Martin2013}. Our concentration argument makes this interpretation explicit by showing that the fixed-data ceiling minimizes parameter volume at each level. Since this ceiling need not be sampling valid, Liu and Martin's recommendation to keep the adjustment for validity small motivates our proposed constructions.

Under the invariant-model assumptions of \citet[Theorem~4]{martin2023fiducial}, the fiducial distribution agrees with the right-Haar posterior and is a \emph{maximal member} of the IM's credal set, that is, a saturating inner probability. \citet[Theorems~2 and~3]{reIM} further relates saturation to probability matching for comonotone rankings and expresses the conditional invariant-model contour in the density-rank form of \eqref{eq:fid-ceiling} when $\nu$ is right Haar measure. In our setting, the same fiducial law saturates every member of the orbit generated by the fixed association.

\citet{hannig_demyst2025} interpret IMs as confidence curves derived from fiducial distributions. Indeed, under the fully identified setting we are working in, many of the identifications made in this paper can be viewed as restatements of their connections from a possibilistic IM perspective on GFD, rather than as a GFD perspective on IMs. Their Theorem~3.3 identifies nested \emph{principle assertions} $A_{\beta,x}$ that determine the IM's belief function. When fiducial inversion requires no conditioning and auxiliary calibration is exact, they establish the existence of a version of the GFD satisfying \[\lPi_x(A_{\beta,x})=\mathrm R_x(A_{\beta,x})=\beta\] at every attained belief level. Their Theorem~4.2 gives at least nominal coverage and, under full compatibility, identifies the actual coverage with the corresponding fiducial content. In the present setting these assertions agree with the cuts $C_{1-\beta}^\gamma(x)$ up to level boundaries that are null under both calibrations. The associated confidence curve is $cc_x^\gamma=1-\pi_x^\gamma$. Thus their fiducial confidence-curve interpretation identifies precisely the agreement on nested assertions expressed here by saturation.

The efficiency-ceiling interpretation has a further direct connection with their discussion. Their Theorem~3.1 places fiducial probability between IM belief and plausibility when the auxiliary draws are almost surely compatible and the predictive random set is calibrated. For a nonatomic auxiliary law, their Theorem~3.2 constructs a nested predictive random set and a fiducial version attaining equality with belief for one fixed observation and assertion. They consequently interpret fiducial probability as a natural bound on IM support, with the additional imprecision in the IM permitting frequentist validity and protection against false confidence \citep{Martin2019FalseConfidence}. In our setting, each $\alpha$-cut of a fixed-data auxiliary-compatible contour has belief at least $1-\alpha$ and hence at least that much fiducial probability. Proposition~\ref{prop:fid-ceiling} connects this bound to efficiency by identifying the density ordering whose cuts attain it with minimum $\nu$-volume under the assumptions of Theorem~\ref{thm:fix-data-ceiling}.

The confidence-curve interpretation emphasizes which fiducial credible regions are also confidence regions. The fiducial orbit perspective provides a partial explanation of the seeming bountifulness of these credible regions. The possibilistic interpretation, however, places emphasis on the calibrated IM as an uncertainty model, with the GFD as a saturating inner. Cui and Hannig themselves note the connection with possibility theory; here that connection is central to how the \emph{inferential output is used}. Committing to one contour fixes a rule for evaluating evidence while leaving open which assertions are considered after observing the data \citep[Corollary~2]{CellaMartin2023Probing}. It provides calibrated lower and upper probabilities for general assertions and validity-preserving uncertainty propagation through its associated calculus \citep{Martin2023b}. Saturation explains where fiducial probability agrees with this reasoning; it does not transfer the IM's general calibration guarantee to unrestricted use of the fiducial probability calculus. The fiducial orbit thus supports exact or approximate calibration of particular fiducial assessments, whereas commitment to one member supplies the richer framework for calibrated uncertainty reasoning. 

\section{Examples}\label{sec:examples}

\subsection{A two-component measurement equation}
\label{sec:metrology-example}

\citet[Remark~2]{hannig2009generalized} emphasizes that generalized fiducial inference is especially natural
when the data-generating mechanism is supplied by the scientific model itself.
In particular, metrology provides a canonical setting in which an unknown true
quantity is related to an observed measurement through a physical measurement
equation, which can then serve directly as the structural equation for
fiducial inversion. Motivated by this connection,
we consider a classical two-component measurement-error model combining an
absolute error floor with multiplicative measurement variation.

A representative example is the analytical-chemistry model of
\citet[equation~(1.4)]{Rocke}, who model a measured concentration through an additive Gaussian
error and a lognormal multiplicative error. We use the same sampling model,
written here in a mean-calibrated parameterization,
\begin{equation}
    \label{eq:metrology-model}
    X
    =
    \theta G+\sigma_0 Z,
    \qquad
    G
    =
    \exp\!\left\{\tau V-\frac{\tau^2}{2}\right\},
    \qquad
    Z,V\stackrel{\mathrm{ind}}{\sim}N(0,1).
\end{equation}
Indeed, the Rocke--Lorenzato equation
$X=\mu\exp(\tau V)+\sigma_0Z$ is obtained by setting
$\mu=\theta\exp(-\tau^2/2)$. Thus the centering in
\eqref{eq:metrology-model} is only a reparameterization, chosen so that
$\theta$ is the mean signal. Whereas the original analytical-chemistry
interpretation has a non-negative concentration parameter, we allow
$\theta\in\mathbb R$, so that $\theta$ may more generally represent a signed
signal or deviation from a target value.

The additive component $\sigma_0 Z$ represents an absolute measurement-error
floor, while the positive random gain $G$ represents multiplicative variation
in measurement sensitivity. Since $E(G)=1$,
\[
    E_\theta(X)=\theta,
    \qquad
    \operatorname{Var}_\theta(X)
    =
    \sigma_0^2+\lambda^2\theta^2,
    \qquad
    \lambda
    :=
    \sqrt{e^{\tau^2}-1}.
\]
Hence $\lambda$ is the coefficient of variation of the random gain.
The additive and multiplicative contributions to the variance are equal when
$\lambda|\theta|=\sigma_0$; additive error dominates below this scale, whereas
relative error dominates for larger signals.

The physical equation \eqref{eq:metrology-model} contains two auxiliary
variables, whereas our theory assumes a fully identified association of the same dimension as the parameter. An equivalent scalar association is obtained from the sampling distribution function $F_\vartheta(x)=E_G\!\left[\Phi\!\left(\frac{x-\vartheta G}{\sigma_0}\right)\right]$ as
\begin{equation}
    \label{eq:metrology-association}
    X
    =
    F_\theta^{-1}\{\Phi(U)\},
    \qquad
    U\sim N(0,1),
\end{equation}
where the inverse is with respect to the data argument. Since
$F_\vartheta(x)$ is strictly decreasing from $1$ to $0$ as $\vartheta$
ranges from $-\infty$ to $\infty$, the compatible auxiliary value is uniquely
defined by
\begin{equation}
    \label{eq:metrology-compatible-u}
    u_{x,\vartheta}
    =
    \Phi^{-1}\{F_\vartheta(x)\}.
\end{equation}
Thus \eqref{eq:metrology-association} satisfies the full-identification and
compatibility assumptions used in the preceding sections. Moreover, the
fiducial law generated by this scalar association agrees with the law obtained
by directly solving the physical measurement equation for $\theta$ with a
fresh pair of measurement errors; see Appendix~\ref{app:metrology-details}.

The sampling density is $p_\vartheta(x)=E_G\!\left[\frac{1}{\sigma_0}\phi\!\left(\frac{x-\vartheta G}{\sigma_0}\right)\right]$.
Take Euclidean parameter volume, $\nu(d\vartheta)=d\vartheta$. If
$\mathrm R_x$ denotes the generalized fiducial distribution induced by
\eqref{eq:metrology-association}, then its density has the particularly simple
form
\begin{equation}
    \label{eq:metrology-fid-density}
    r_x(\vartheta)
    =
    p_{\kappa\vartheta}(x),
    \qquad
    \kappa
    :=
    e^{\tau^2}
    =
    1+\lambda^2.
\end{equation}
Consequently, the association-induced correction in
Proposition~\ref{prop:fid-ceiling} is
\begin{equation}
    \label{eq:metrology-J}
    J(x,\vartheta)
    =
    \frac{p_{\kappa\vartheta}(x)}
         {p_\vartheta(x)}
    =
    E_\vartheta(G\mid X=x).
\end{equation}
The fiducial adjustment therefore has a direct measurement interpretation:
the likelihood is corrected by the conditional measurement sensitivity under
the candidate value. The identities
\eqref{eq:metrology-fid-density}--\eqref{eq:metrology-J} are derived in
Appendix~\ref{app:metrology-details}.

\paragraph{Contours under comparison.}
We compare the constructions developed above, specialized to the present
association.  For the standard-normal auxiliary,
maximum specificity gives
\begin{equation}
    \label{eq:metrology-ms}
    \pi_x^{\mathrm{ms}}(\vartheta)
    =
    2\min\{F_\vartheta(x),1-F_\vartheta(x)\}.
\end{equation}
Its cuts are the equal-tailed confidence sets obtained by inversion of the
sampling distribution; after the reparameterization above, their endpoints solve the same tail equations as the exact single-measurement intervals of \citet[equations~(3.1)--(3.2)]{Rocke} when the two error scales are known.

By Proposition~\ref{prop:fid-ceiling} and
\eqref{eq:metrology-fid-density}, the fixed-data ceiling is
\begin{equation}
    \label{eq:metrology-ceiling}
    \pi_x^{\mathrm{ceil}}(\vartheta)
    =
    \mathrm R_x
    \{p_{\kappa\Theta}(x)\le p_{\kappa\vartheta}(x)\}.
\end{equation}
Let $M(x)=\sup_{\eta\in\mathbb R}p_\eta(x)$.
Since $\sup_\eta p_{\kappa\eta}(x)=M(x)$, the normalized likelihood and
fiducial-likelihood scores become
\begin{equation}
    \label{eq:metrology-relative-scores}
    R_x^{\mathrm{LR}}(\vartheta)
    =
    \frac{p_\vartheta(x)}{M(x)},
    \qquad
    RF_x(\vartheta)
    =
    \frac{p_{\kappa\vartheta}(x)}{M(x)}.
\end{equation}
Their candidatewise validifications give, respectively, the LR contour
$\pi_x^{\mathrm{LR}}$ of Section~\ref{sec:validified-LR} and the
fiducial LR/localized contour $\pi_x^{\mathrm{loc}}$ of
Section~\ref{sec:loc-repair}. We also include the directly validified
ceiling $\pi_x^{\mathrm{Vc}}$ from Section~\ref{sec:direct-rep}, and finally, the $w$-averaged optimal contour $\pi^\star_{x,w}$ in
Theorem~\ref{thm:prs-w-choquet-ceiling}, using  $w=\operatorname{Unif}\left[0,\frac{4\sigma_0}{\lambda}\right]$.
This range of $w$ spans the transition from additive- to multiplicative-error
dominance. 

\paragraph{Concentration comparison.}
For every method $m$, we evaluate the same integrated-contour loss $    L_m(x)
    =
    \int_{\mathbb R}
        \pi_x^m(\vartheta)\,d\vartheta$
and its pointwise sampling risk $R_m(\theta)=E_\theta\{L_m(X)\}$.
 We report $RR_m(\theta)=R_m(\theta)/
R_{\mathrm{ceil}}(\theta)$.

\begin{figure}[ht!]
    \centering
    \begin{subfigure}[t]{0.49\textwidth}
        \centering
        \includegraphics[width=\linewidth]{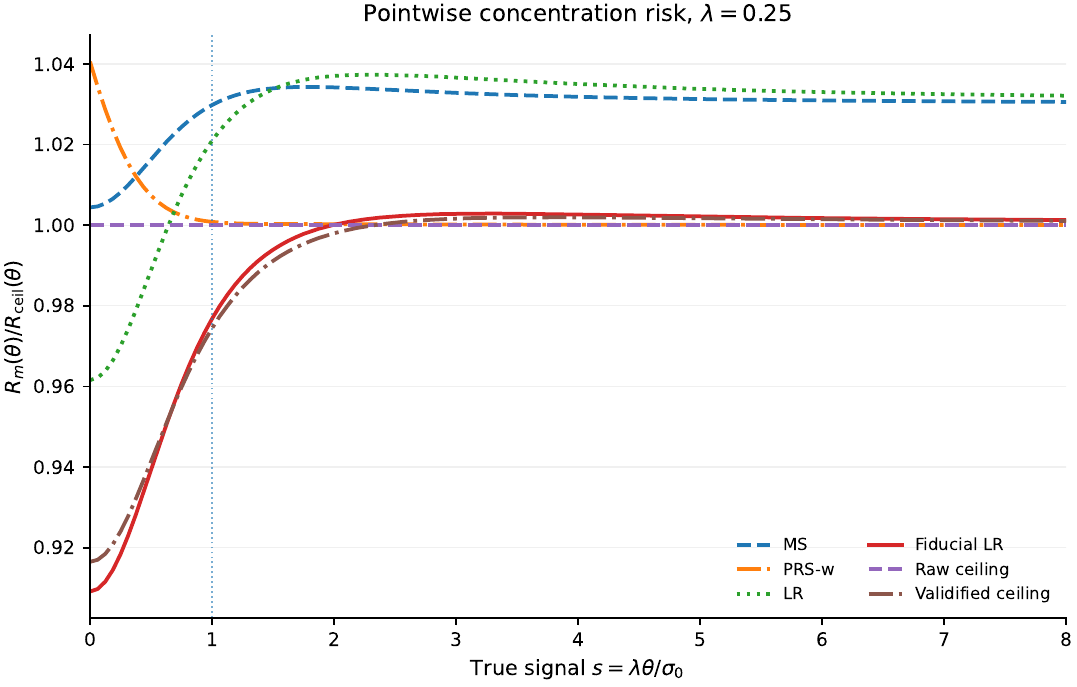}
        \caption{$\lambda=0.25$.}
        \label{fig:metrology-risk-025}
    \end{subfigure}
    \hfill
    \begin{subfigure}[t]{0.49\textwidth}
        \centering
        \includegraphics[width=\linewidth]{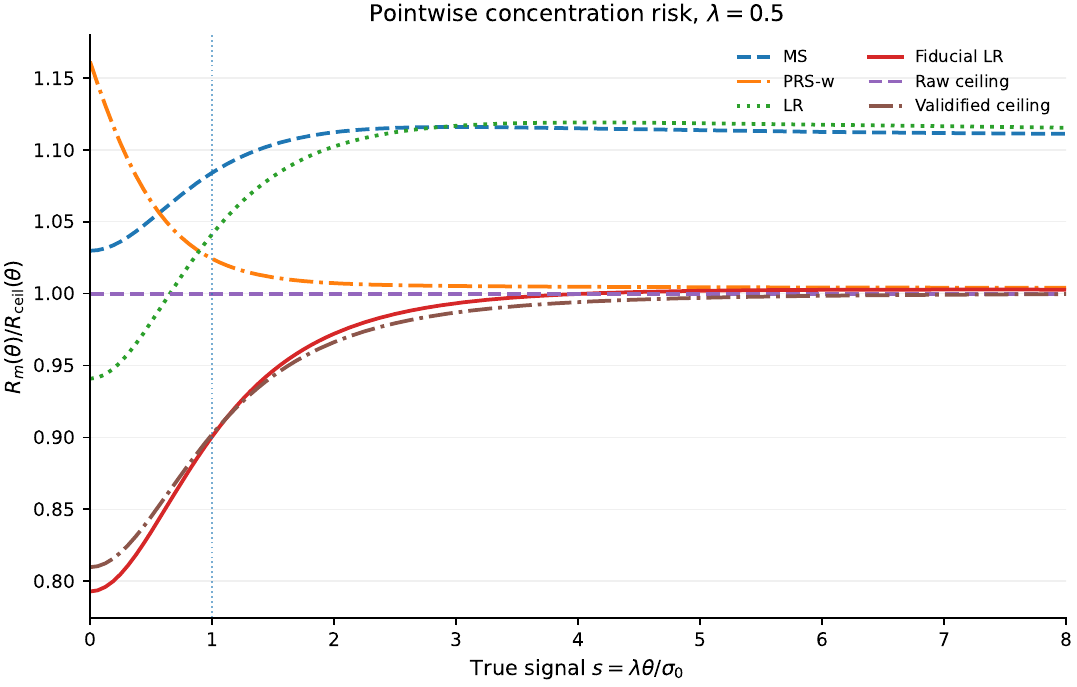}
        \caption{$\lambda=0.5$.}
        \label{fig:metrology-risk-05}
    \end{subfigure}
    \hfill
    \begin{subfigure}[t]{0.49\textwidth}
        \centering
        \includegraphics[width=\linewidth]{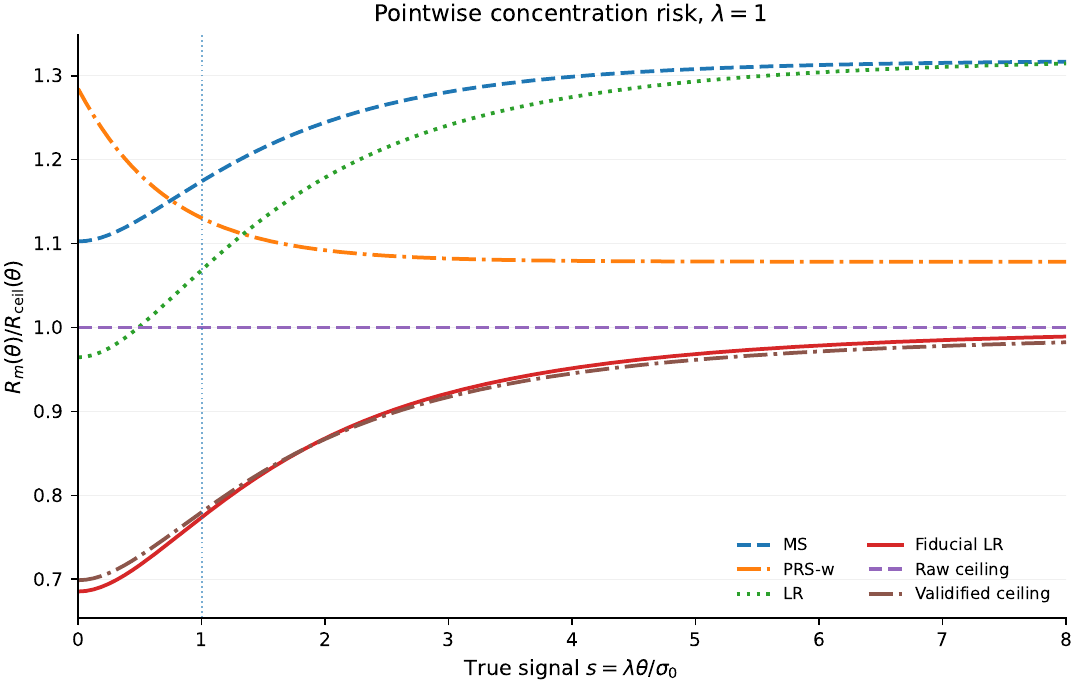}
        \caption{$\lambda=1$.}
        \label{fig:metrology-risk-1}
    \end{subfigure}
    \caption{
        Pointwise concentration risk in the two-component measurement model,
        relative to the raw fixed-data ceiling, for $\sigma_0=1$ and
        increasing gain variation $\lambda$. Lower values indicate
        greater contour concentration. 
    }
    \label{fig:metrology-risks}
\end{figure}

We set $\sigma_0=1$, consider $\lambda\in\{0.25,0.5,1\}$, and index the
true parameter by $s=\frac{\lambda\theta}{\sigma_0}$.
Hence $|s|=1$ marks the point at which the additive and multiplicative
variance components are equal. The risks in
Figure~\ref{fig:metrology-risks} are evaluated by Monte Carlo under the
model \eqref{eq:metrology-model}; the contour integrals and the
candidatewise calibration distributions are evaluated numerically. Common
random numbers are used across the competing procedures. Further numerical
details are given in Appendix~\ref{app:metrology-details}.

Figure~\ref{fig:metrology-risks} shows a common pattern across the three
gain-variation levels. MS and LR attain similar relative risks at larger
signals, but LR is more concentrated near zero and falls below the raw
ceiling. The $w$-optimal contour becomes close to
the ceiling in relative risk at larger signals for $\lambda=0.25$ and
$0.5$, whereas a visible gap remains over the displayed high-signal range
for $\lambda=1$. This high signal gap reflects the different efficiency priorities encoded by the chosen design, and which seem to depend on $\lambda$.
Fiducial LR and the directly validified ceiling have closely tracking
risk curves, with fiducial LR slightly more efficient near zero and
direct validification generally more efficient thereafter. Both improve
substantially on the ceiling near the additive-error floor and remain
below or close to it as the signal increases, illustrating that
sampling-valid repair need not incur a concentration penalty.
Nevertheless, for $\lambda=0.25$, the $w$-optimal contour slightly
outperforms both repairs at larger signals, so neither repair is
uniformly preferred across the comparisons. Overall, increasing
$\lambda$ magnifies both the large-signal separation of MS and ordinary
LR from the parameter-volume-based constructions and the near-zero
gains from the fiducial repairs. 

\begin{figure}[ht!]
    \centering
    \includegraphics[width=0.49\textwidth]{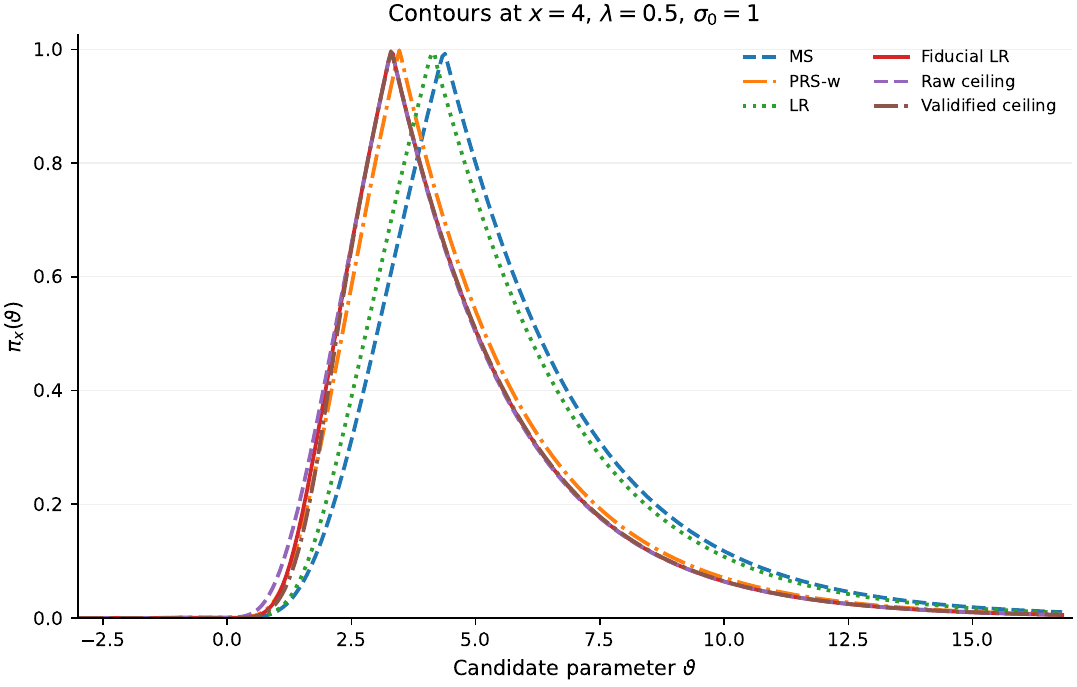}
    \caption{
        Contours at $x=4$ for $\sigma_0=1$ and $\lambda=0.5$.
        The $w$-optimal PRS uses
        $w=\operatorname{Unif}[0,8]$.
    }
    \label{fig:metrology-contours}
\end{figure}

Figure~\ref{fig:metrology-contours} shows the corresponding contours for
$x=4$, $\sigma_0=1$, and $\lambda=0.5$. This observation lies well beyond the
additive-error floor and makes the difference between likelihood and
fiducial-volume ordering visible. From
\eqref{eq:metrology-relative-scores}, $\widehat\theta_{\mathrm F}(x)=\widehat\theta_{\mathrm L}(x)/\kappa$,
so the fiducial score is maximized at a smaller parameter value than the
ordinary likelihood whenever $\lambda>0$. Numerically, at $x=4$ the two
maxima are approximately $3.30$ and $4.13$, respectively. The fiducial LR and
directly validified ceiling remain close, as reflected in their
similar risks.

\subsection{An exponential scale model}
\label{sec:exponential-example}

We return to the exponential scale model of
Section~\ref{sec:predictive-contour-choice} to compare the resulting
contours under ordinary parameter length. We take $n\geq 3$ so that
the concentration risks compared below are finite.

\paragraph{Contours under comparison.}
For $k\in\{n-1,n,n+1\}$, write
\[
    h_k(u)=u^k e^{-u},
    \qquad
    \gamma_{n,k}(u)=\mathrm P_U\{h_k(U)\leq h_k(u)\},
    \qquad U\sim\operatorname{Gamma}(n,1).
\]
Auxiliary maximum specificity uses $h_{n-1}$, while the ceiling and
the design-averaged construction both use $h_{n+1}$ by
Corollary~\ref{cor:factoriz}. To identify the likelihood-based
contours, observe that the relative likelihood and relative fiducial
likelihood are
\[
    R_x^{\mathrm{LR}}(\vartheta)
      =\frac{h_n(x/\vartheta)}{h_n(n)},
    \qquad
    RF_x(\vartheta)
      =\frac{h_{n+1}(x/\vartheta)}{h_{n+1}(n+1)}.
\]
Under $\mathrm P_\vartheta$, the pivot $X/\vartheta$ has the same
gamma law used to define $\gamma_{n,k}$. Since the denominators above
are constant, candidatewise validification amounts to calibrating each
numerator under this gamma law. In particular, the localized fiducial
LR contour coincides with the ceiling. Collecting the six constructions
gives

\begin{equation}
    \label{eq:exponential-contours}
    \begin{aligned}
        \pi_x^{\mathrm{ms}}(\vartheta)
        &=\gamma_{n,n-1}(x/\vartheta),\\
        \pi_x^{\mathrm{LR}}(\vartheta)
        &=\gamma_{n,n}(x/\vartheta),\\
        \pi_x^{\mathrm{ceil}}(\vartheta)
        =\pi^\star_{x,w}(\vartheta)
        =\pi_x^{\mathrm{loc}}(\vartheta)
        =\pi_x^{\mathrm{Vc}}(\vartheta)
        &=\gamma_{n,n+1}(x/\vartheta).
    \end{aligned}
\end{equation}
Indeed, each $h_k(U)$ has a continuous distribution, so
$\gamma_{n,k}(U)\sim\operatorname{Unif}(0,1)$.
All three contours are normalized and exactly sampling-calibrated;
in particular, direct validification leaves the ceiling unchanged.

\paragraph{Concentration comparison.}
For the contour generated by $\gamma_{n,k}$, changing variables from
$\vartheta$ to $u=x/\vartheta$ gives
\begin{equation}
    \label{eq:exponential-risk}
    \begin{aligned}
        L_{n,k}(x)
        &=
        \int_0^\infty
            \gamma_{n,k}(x/\vartheta)\,d\vartheta
        =x\,c_{n,k},\\
        R_{n,k}(\theta)
        &=n\theta\,c_{n,k},
        \qquad
        c_{n,k}
        :=
        \int_0^\infty\frac{\gamma_{n,k}(u)}{u^2}\,du.
    \end{aligned}
\end{equation}
Unlike the measurement example, all competing contours here belong
to the single-PRS class. Theorem~\ref{thm:fix-data-ceiling} therefore
implies that the ceiling minimizes cut length at every level and
observation, and hence concentration risk at every $\theta>0$ within
this class. Among exact scale-equivariant confidence sets based on $X$, the ceiling has shortest Euclidean parameter length. The pivot $2X/\theta\sim\chi^2_{2n}$ gives the comparison studied by \citet[Example~6.12.1]{lehmann2005testing} for a normal variance. Inverting the two-sided uniformly most powerful unbiased (UMPU) tests gives the LR ordering $u^ne^{-u}$ \citep[Section~4.2]{lehmann2005testing}. For a pivot acceptance interval $[a,b]$ on the $U$ scale, unbiasedness gives $a^ne^{-a}=b^ne^{-b}$, whereas minimum Euclidean parameter length gives the ordering $u^{n+1}e^{-u}$.

\begin{figure}[ht!]
    \centering
    \begin{subfigure}[t]{0.49\textwidth}
        \centering
        \includegraphics[width=\linewidth]{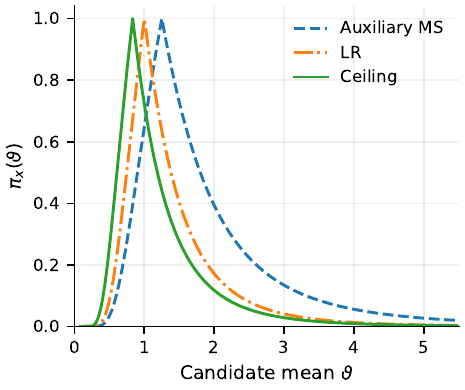}
        \caption{Contours for $n=5$ and $x=5$.}
        \label{fig:exponential-contours}
    \end{subfigure}
    \hfill
    \begin{subfigure}[t]{0.49\textwidth}
        \centering
        \includegraphics[width=\linewidth]{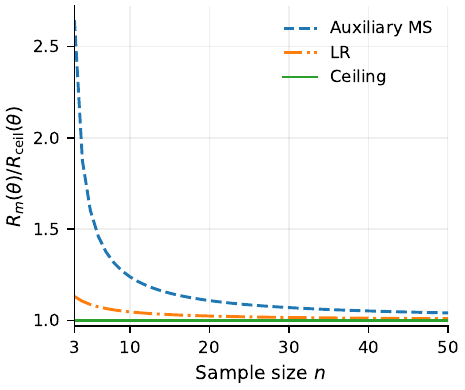}
        \caption{Relative concentration risk.}
        \label{fig:exponential-risk}
    \end{subfigure}
    \caption{
        Exponential scale model with Euclidean parameter volume.
        The ceiling, $w$-optimal PRS, fiducial LR, and directly
        validified ceiling coincide. 
    }
    \label{fig:exponential-comparison}
\end{figure}

Relative risks against the ceiling are independent of $\theta$, and depend only on sample size $n$.
For $n=5$, the relative risks are $RR_{\mathrm{ms}}\approx1.6023$ and $RR_{\mathrm{LR}}\approx1.0873$. Figure~\ref{fig:exponential-contours} shows the three contours, whose
maximizers are $x/(n-1)$, $x/n$, and $x/(n+1)$, respectively.
Figure~\ref{fig:exponential-risk} shows the diminishing risk differences
over the displayed sample sizes.

The comparison depends on the chosen volume measure. With logarithmic/Haar
volume, $\nu(d\vartheta)=d\vartheta/\vartheta$, the pullback cost is
instead $K_x^\nu(u)=1/u$, so $f_U(u)/K_x^\nu(u)\propto h_n(u)$ and
ordinary LR coincides with the corresponding ceiling. Thus the
Euclidean improvement reflects alignment with the specified
concentration loss, rather than a parameterization-independent
preference.

\section{Further discussions}\label{sec:discussions}
The preceding constructions raise broader questions about efficiency criteria and how to choose among contours. We have seen potential efficiency gain by leaving the single predictive contour paradigm; might that be achievable also when the designed-averaged risk is of actual interest, rather than just a device for construction? That is, can we keep the design-averaged risk criteria, but permit candidate-local contour constructions?  
\subsection{An unrestricted averaged-risk benchmark and renormalization}
Theorem~\ref{thm:prs-w-choquet-ceiling} optimizes (averaged) concentration within the single predictive contour class. A broader benchmark can be obtained by removing this restriction and minimizing the same averaged risk over all measurable sampling-valid surfaces. Retaining the design measure $w$ and parameter-volume measure $\nu$, write \[m_w(x)=\int_\Theta p_\eta(x)w(d\eta),\qquad R_w(\pi)=\int_\mathcal Xm_w(x)\int_\Theta\pi_x(\vartheta)\nu(d\vartheta)dx.\]
For this discussion, assume that $w$ is finite, the mixture is positive and finite, and the risks under comparison are finite. At each fixed confidence level and for ordinary parameter length, this is the design-averaged expected-length problem studied by \citet{Pratt01091961}. The optimization problem again reduces to Neyman--Pearson: for each candidate and confidence level, sampling validity requires an acceptance region with prescribed $P_\vartheta$-probability, while its contribution to the objective is its $m_w$-weighted sample-space volume. The Neyman--Pearson ordering is therefore $p_\vartheta/m_w$, giving \[\pi_x^{m_w}(\vartheta)=P_\vartheta\left\{ \frac{p_\vartheta(X)}{m_w(X)}\le \frac{p_\vartheta(x)}{m_w(x)}\right\}.\] When the candidatewise score distributions are continuous, this surface minimizes the design-risk $R_w$ over all sampling-valid surfaces.
Its cuts minimize averaged expected $\nu$-volume separately at every level and are nested because they threshold a common score for each candidate.

The concentration criterion integrates these fixed-level objectives over confidence levels for the specified volume measure $\nu$. The mixture denominator thus arises from the efficiency criterion, without requiring a prior-to-posterior interpretation. As before, the measure $w$ specifies where efficiency is prioritized, while validity holds throughout the parameter space.

The unrestricted optimal surface, whilst valid (i.e. a p-value surface), it need not be normalized, and hence not a possibility contour. Its higher cuts may be empty, and the integrated-volume functional then lacks its interpretation as an upper expectation for a normalized possibility measure. This motivates a sort of \emph{possibilistic projection}: minimize the same risk, but over a restricted class of valid and normalized contours. We call it ``projection" here because the unrestricted optimal solution inspires the functional form of the possibilistic class.

Without further restrictions, this problem admits degenerate solutions. If $\nu$ is non-atomic, fixing any $\vartheta_0$ and replacing the unrestricted optimum by \[\tilde \pi_x(\vartheta)=\max\{\pi_x^{m_w}(\vartheta),1_{\{\vartheta_0\}}(\vartheta)\}\] preserves validity and enforces normalization without changing $R_w$. The volume criterion cannot detect the artificial spike. Consequently, normalization alone does not identify a useful repair: the class of admissible construction must also be specified.

A natural restriction retains the ratio-score form of the unrestricted optimum. Fix a positive score family $g_x(\vartheta)$ and consider $\pi^q=\mathcal V[g/q]$, where the denominator $q(x)>0$ is shared across candidates, and validification remains under $P_\vartheta$. Taking $g_x(\vartheta)=p_\vartheta(x)$ recovers the likelihood-based class and functional form of the unrestricted benchmark; other fixed scores permit association-induced geometry to enter the construction, such as the fiducial likelihood.

The question is now which denominators produce a normalized contour. To investigate this, consider the upper undpoint of the score distribution under each candidate: \[M_q(\vartheta)=\esssup_{X\sim P_\vartheta} \frac{g_X(\vartheta)}{q(X)}, \qquad b_q(\vartheta)=M_q(\vartheta)^{-1}.\] Assume that these endpoints are finite and nonzero, that the sampling laws have common full support, and that $g_\cdot(\vartheta)$ and $q$ are positive and continuous there. The essential upper bound then holds at every observation: $b_q(\vartheta)g_x(\vartheta)\le q(x)$. Multiplication by the positive candidate-specific constant $b_q(\vartheta)$ does not change the probability comparison defining $\pi^q$. Moreover, the calibrated values reaches one precisely when the observed score reaches its upper endpoint. Hence \[\pi_x^q(\vartheta)=1 \Longleftrightarrow b_q(\vartheta)g_x(\vartheta)=q(x).\] Requiring a fully plausible candidate at every observation therefore forces \[q(x)=\max_{\vartheta\in\Theta} b_q(\vartheta)g_x(\vartheta).\] Thus $q$ is a \emph{weighted score envelope}.

Conversely, for a positive weight $b$, define \[T_x^b=\frac{b(\vartheta)g_x(\vartheta)}{\sup_{\eta \in \Theta}b(\eta)g_x(\eta)}, \qquad \pi^b=\mathcal V [T^b].\] A finite, attained envelope makes the score bounded by one and equal to one at a maximizer, sot its validification is valid and normalized. Since $B(\vartheta)$ cancels from the candidate probability comparison, weighted envelopes charactarize attained normalization within the stipulated regular common-denominator class. With $b\equiv1$, choosing $g_x(\vartheta)=p_\vartheta(x)$ gives the ordinary likelihood-ratio contour, while choosing the fiducial score $g_x(\vartheta)=s_x(\vartheta)$ gives the localized, fiducial likelihood-ratio contour. 

For a specified admissible family $\mathcal B_g$ of positive weights, the possibilistic projection problem consequently becomes \[\inf_{b\in\mathcal B_g} R_w(\pi^b),\] subject to the required regularity, finiteness, and attainment conditions. 

 \subsection{Partial priors and efficiency}
 For the likelihood-score family, a bounded positive weight $b$ can be rescaled to satisfy $\sup_\vartheta b(\vartheta)=1$, giving it the mathematical form of a prior possibility contour. The resulting normalized score is the likelihood-times-prior score used by \citet{Martin2023b}. His Remark 4 motivates supremum normalization by its ability to preserve the required possibility structure; the argument above supplies a complementary necessity result within the specified common-denominator class. More specifically, our candidatewise validification agrees with his Section 6.4 construction, rather than the general construction that calibrates under an upper joint distribution incorporating the partial prior.

 Here, the weight $b$ may instead serve purely as a construction parameter. Its numerator factor cancels during candidatewise validification, so it affects the final contour only through the common envelope and the resulting sampling orderings. Prior-shaped weighting thus offers a mechanism for potential efficiency improvement without requiring prior belief. Any gain should be assessed through the concentration-risk function, or through an explicitly intended averaged criterion.

 Existence and characterization of optimal weights remain further questions. Unlike the unrestricted Neyman--Pearson problem, changing the envelope also changes the sampling distributions and hence couples the candidate orderings, so the problem no longer separates into independent Neyman-Pearson optimizations. The formulation here nevertheless places the likelihood-ratio and fiducial likelihood-ratio constructions within a common efficiency program and identifies how their weighting might be improved.
 
\section{Concluding remarks}

The constructions developed here connect parameter-space concentration, fiducial inference, and the choice of valid IMs. The fiducial density emerges as the score solving the fixed-data volume problem. Moreover, exact single predictive contours share a common fiducial law as a saturating inner probability, while generating different geometries of calibrated confidence regions. Design averaging selects an optimal member of this class for the specified averaged concentration criterion. Localization uses the same parameter-volume geometry through candidate-specific auxiliary contours to obtain further valid constructions. These connections give an efficiency-oriented interpretation of familiar fiducial ingredients and clarify how they can guide IM construction.

The explicit optimality results rely on the smooth, fully identified association considered here. Beyond this setting, inversion may produce sets of solutions, and observations may restrict the compatible auxiliary values. The random-set formulation accommodates nonunique inversion, and validity can be retained under suitable calibration and nonemptiness conditions \citep{Martin2013}, although exactness may be lost. The efficiency analysis is more delicate: overlap among parameter images can prevent their volume from being expressed through the additive Jacobian cost used here. The common saturating fiducial law is likewise no longer automatic. Appropriate dimension reduction or conditioning may recover the present structure in some models \citep{MartinConditional2015}. Identifying analogous efficiency benchmarks more generally remains an open direction.

Our efficiency conclusions also concern a particular loss and parameter-volume measure. Concentration assesses the extent of the reported uncertainty, while other Choquet losses can assess its proximity to the true parameter \citep{raner2026efficiency}. Extending the ceiling and design-averaged constructions to such criteria would broaden their applicability, but may require different optimization principles.

The fiducial connections emphasize the continuity between these IM constructions and established approaches to inference. Familiar probability distributions and confidence regions remain part of the resulting uncertainty description. The additional imprecision permits a suitably calibrated IM to provide lower and upper probabilities for general assertions, with frequentist control of misleadingly strong support even for assertions selected after seeing the data \citep{CellaMartin2023Probing}. Its possibility calculus also preserves validity when uncertainty is propagated to functions of the parameter, including through marginalization \citep{Martin2023b}. Agreement with fiducial probability on distinguished nested regions thus fits within a broader framework for uncertainty reasoning. The price is giving up a single additive probability assessment for every assertion; the gain is a calibrated calculus that extends the use of familiar inferential constructions beyond those regions.

\section*{Acknowledgments}
I thank Ryan Martin for helpful comments on an earlier draft.

\section*{Declaration of generative AI use}
I used ChatGPT (OpenAI) to help with writing and editing, to identify and review relevant literature, and to test ideas. I take responsibility for the content of this paper.

\bibliographystyle{plainnat}
\bibliography{sample}

\appendix
\section{Details for the two-component measurement example}
\label{app:metrology-details}

We collect the calculations used in
Section~\ref{sec:metrology-example}. Conditional on the random gain $G$,
\[
    X\mid G
    \sim
    N(\theta G,\sigma_0^2),
\]
so
\[
    F_\vartheta(x)
    =
    E_G\!\left[
        \Phi\!\left(
            \frac{x-\vartheta G}{\sigma_0}
        \right)
    \right],
    \qquad
    p_\vartheta(x)
    =
    E_G\!\left[
        \frac{1}{\sigma_0}
        \phi\!\left(
            \frac{x-\vartheta G}{\sigma_0}
        \right)
    \right].
\]
Differentiation under the expectation gives
\begin{equation}
    \label{eq:app-r-density}
    -\partial_\vartheta F_\vartheta(x)
    =
    E_G\!\left[
        \frac{G}{\sigma_0}
        \phi\!\left(
            \frac{x-\vartheta G}{\sigma_0}
        \right)
    \right]
    >0.
\end{equation}
Together with
$F_\vartheta(x)\to1$ as $\vartheta\to-\infty$ and
$F_\vartheta(x)\to0$ as $\vartheta\to\infty$, this proves global
invertibility of the scalar association
\eqref{eq:metrology-association}.

If $\Theta_x=\theta_x(U)$ for $U\sim N(0,1)$, then, since
$\vartheta\mapsto F_\vartheta(x)$ is decreasing,
\[
\begin{aligned}
    \mathrm R_x((-\infty,\vartheta])
    &=
    P\{\Theta_x\leq\vartheta\}  \\
    &=
    1-F_\vartheta(x).
\end{aligned}
\]
The same law follows directly from the physical measurement equation.
Indeed, for independent copies $(Z',G')$ of the measurement errors,
\[
    \widetilde\Theta_x
    =
    \frac{x-\sigma_0Z'}{G'}
\]
satisfies
\[
\begin{aligned}
    P(\widetilde\Theta_x\leq\vartheta)
    &=
    E_G
    \left[
        \Phi\!\left(
            \frac{\vartheta G-x}{\sigma_0}
        \right)
    \right] \\
    &=
    1-F_\vartheta(x).
\end{aligned}
\]
Thus the scalar reduction preserves the inverse-equation fiducial law.

To obtain its density, write
\[
    g(v)
    =
    \exp\!\left\{\tau v-\frac{\tau^2}{2}\right\},
    \qquad
    \kappa=e^{\tau^2}.
\]
Since
\[
    g(v)\phi(v)=\phi(v-\tau),
\]
a change of variable gives, for any integrable $h$,
\begin{equation}
    \label{eq:lognormal-tilt}
    E\{Gh(G)\}
    =
    E\{h(\kappa G)\}.
\end{equation}
Applying \eqref{eq:lognormal-tilt} to
\eqref{eq:app-r-density} gives
\[
\begin{aligned}
    r_x(\vartheta)
    &=
    -\partial_\vartheta F_\vartheta(x)\\
    &=
    E_G\!\left[
        \frac{G}{\sigma_0}
        \phi\!\left(
            \frac{x-\vartheta G}{\sigma_0}
        \right)
    \right] \\
    &=
    p_{\kappa\vartheta}(x),
\end{aligned}
\]
which proves \eqref{eq:metrology-fid-density}.

The same calculation gives the fiducial Jacobian interpretation. From
\[
    F_\vartheta\{a(\vartheta,u)\}
    =
    \Phi(u),
\]
implicit differentiation with respect to $\vartheta$ yields
\[
    \partial_\vartheta a(\vartheta,u)
    =
    \frac{
        r_{a(\vartheta,u)}(\vartheta)
    }{
        p_\vartheta\{a(\vartheta,u)\}
    }.
\]
Consequently,
\[
    J(x,\vartheta)
    =
    \frac{r_x(\vartheta)}
         {p_\vartheta(x)}
    =
    \frac{p_{\kappa\vartheta}(x)}
         {p_\vartheta(x)}.
\]
On the other hand,
\[
    p_\vartheta(x)
    E_\vartheta(G\mid X=x)
    =
    E_G\!\left[
        \frac{G}{\sigma_0}
        \phi\!\left(
            \frac{x-\vartheta G}{\sigma_0}
        \right)
    \right],
\]
so
\[
    J(x,\vartheta)
    =
    E_\vartheta(G\mid X=x).
\]

For the $w$-optimal PRS, differentiation of
\[
    F_{\theta_x(u)}(x)
    =
    \Phi(u)
\]
gives
\[
    K_x(u)
    =
    |\theta_x'(u)|
    =
    \frac{\phi(u)}
         {r_x\{\theta_x(u)\}}.
\]
Hence
\[
\begin{aligned}
    \bar K_w(u)
    &=
    \phi(u)D_w(u),\\
    D_w(u)
    &:=
    \int
    E_\eta
    \left[
        \frac{1}
        {r_X\{\theta_X(u)\}}
    \right]
    w(d\eta),
\end{aligned}
\]
and therefore
\[
    \rho_w(u)
    =
    \frac{1}{D_w(u)}.
\]
This form is used for the numerical evaluation of the $w$-optimal
single-PRS contour.

The ordinary maximum-specificity contour follows directly from the standard
normal auxiliary density:
\[
\begin{aligned}
    \gamma^{\mathrm{ms}}(u)
    &=
    P\{\phi(U)\leq\phi(u)\}\\
    &=
    P(|U|\geq|u|)\\
    &=
    2\Phi(-|u|).
\end{aligned}
\]
Substituting \eqref{eq:metrology-compatible-u} yields
\[
    \pi_x^{\mathrm{ms}}(\vartheta)
    =
    2\min\{F_\vartheta(x),1-F_\vartheta(x)\}.
\]

For the numerical comparison, the sampling model is standardized to
$\sigma_0=1$. The concentration loss
\[
    L_m(x)
    =
    \int_{\mathbb R}\pi_x^m(\vartheta)\,d\vartheta
\]
is evaluated numerically for each method and interpolated as a function of
$x$. The outer expectation in
$R_m(\theta)=E_\theta\{L_m(X)\}$ is then estimated by paired Monte Carlo
simulation from \eqref{eq:metrology-model}, using common random numbers across
methods. The calculation used $2^{20}$ outer Monte Carlo draws at each
displayed true value; because the contour losses and calibration functions
were precomputed numerically, these draws require only evaluation and
interpolation of the stored loss functions rather than reconstruction of all
contours.

As two numerical checks, when $\tau=0$ the model reduces to
$X=\theta+\sigma_0Z$ and all six constructions coincide with
\[
    \pi_x(\vartheta)
    =
    2\Phi\!\left(
        -\frac{|x-\vartheta|}{\sigma_0}
    \right).
\]
For arbitrary $\tau$, $r_0(\vartheta)$ is symmetric and strictly decreasing
in $|\vartheta|$, so the fixed-data fiducial ceiling and auxiliary
maximum-specificity contour coincide at $x=0$:
\[
    \pi_0^{\mathrm{ceil}}
    =
    \pi_0^{\mathrm{ms}}.
\]
For contrast, the contour comparison in Figure~\ref{fig:metrology-contours} is therefore shown at $x=4$.

\section{Design measures on non-compact parameter spaces}
\label{app:noncompact-design}

The design-averaged construction in
Theorem~\ref{thm:prs-w-choquet-ceiling} assumes that the averaged
volume cost $\bar K_w(u)$ in \eqref{eq:Kbar} is positive and finite
$\mathrm P_U$-almost everywhere and that the risks under comparison
are finite. On a non-compact parameter space, however, a natural
design measure, such as Lebesgue or Haar measure, may have infinite
mass, and the corresponding integral may diverge. Since the
Neyman--Pearson ordering depends only on the relative shape of
$\bar K_w$ over $u$, a meaningful limiting construction may still
be obtained by truncation.

To make this precise, let $C_n\uparrow\Theta$ be an increasing
exhaustion by compact sets with $0<w(C_n)<\infty$, and define the
truncated volume cost and corresponding score by
\[
    \bar K_{w,n}(u)=\int_{C_n}E_\eta[K_X(u)]\,w(d\eta),
    \qquad
    \rho_{w,n}(u)=\frac{f_U(u)}{\bar K_{w,n}(u)}.
\]
For every $n$, this is a finite-design version of the original
optimization problem. Whenever the hypotheses of
Theorem~\ref{thm:prs-w-choquet-ceiling} hold for the restricted
design, the ordering of $\rho_{w,n}$ determines an optimal auxiliary
contour. Since multiplication of $\rho_{w,n}$ by a positive constant
leaves this ordering and the resulting contour unchanged, the
relevant limiting question is whether the auxiliary orderings
stabilize as $C_n\uparrow\Theta$.

The factorization in Corollary~\ref{cor:factoriz} gives the simplest
example. If $K_x(u)=c(x)K_0(u)$, then, for every admissible
truncation $C_n$,
\[
    \bar K_{w,n}(u)
    =K_0(u)\int_{C_n}E_\eta[c(X)]\,w(d\eta)
    =A_nK_0(u),
\]
where $0<A_n<\infty$ does not depend on $u$. Consequently,
$\rho_{w,n}\propto f_U/K_0$ for every $n$. The induced contour is
therefore independent of the truncation, even if $A_n\to\infty$.
In this case the divergence of the design integral occurs only
through an overall scale factor and has no effect on the optimal
auxiliary ordering. More generally, whenever the truncated scores
eventually induce a common ordering along a given exhaustion, the
contours stabilize along that exhaustion even without exact
factorization.

Compact truncation can also be chosen deliberately as part of the
efficiency criterion. In many applications, only a scientifically
relevant region $C\subseteq\Theta$ is of primary interest, even
though the statistical model itself is defined on a much larger
parameter space. Replacing the full design integral by
\[
    \bar K_{w,C}(u)=\int_C E_\eta[K_X(u)]\,w(d\eta)
\]
then asks for a valid auxiliary contour whose concentration is most
efficient, in the design-averaged sense, over the region $C$. The
truncation specifies where inferential efficiency is to be prioritized.

This interpretation should be distinguished from assigning a prior
distribution to $\theta$. Although the design measure enters the
optimization mathematically through an average over parameter values,
it need not represent beliefs about where the true parameter lies.
Restricting the design to $C$ makes no assertion that $\theta\in C$,
nor does it assign any prior probability to that event. The validity
requirement remains the same over the full parameter space. Only the
criterion used to compare valid procedures is localized to $C$.
In this sense, validity remains global while efficiency is designed
locally.

The design measure can therefore be understood as encoding scientific
or experimental priorities. For example, one may choose to optimize
concentration over a range of effect sizes judged scientifically
meaningful, while still requiring the resulting contour to satisfy
its sampling-validity guarantee for all $\theta\in\Theta$. Compact
truncation thus has two complementary roles. It can regularize a
design problem on a non-compact parameter space, or it can be
introduced intentionally to localize efficiency to the part of the
parameter space of greatest scientific relevance.

\end{document}